# FIRST EXIT TIMES FOR LÉVY-DRIVEN DIFFUSIONS WITH EXPONENTIALLY LIGHT JUMPS[1]

By Peter Imkeller, Ilya Pavlyukevich and Torsten Wetzel

*Humboldt—Universität zu Berlin*

We consider a dynamical system described by the differential equation $\dot{Y}_t = -U'(Y_t)$ with a unique stable point at the origin. We perturb the system by the Lévy noise of intensity $\varepsilon$ to obtain the stochastic differential equation $dX_t^\varepsilon = -U'(X_{t-}^\varepsilon)\,dt + \varepsilon\,dL_t$. The process $L$ is a symmetric Lévy process whose jump measure $\nu$ has exponentially light tails, $\nu([u,\infty)) \sim \exp(-u^\alpha)$, $\alpha > 0$, $u \to \infty$. We study the first exit problem for the trajectories of the solutions of the stochastic differential equation from the interval $(-1,1)$. In the small noise limit $\varepsilon \to 0$, the law of the first exit time $\sigma_x$, $x \in (-1,1)$, has exponential tail and the mean value exhibiting an intriguing phase transition at the critical index $\alpha = 1$, namely, $\ln \mathbf{E}\sigma \sim \varepsilon^{-\alpha}$ for $0 < \alpha < 1$, whereas $\ln \mathbf{E}\sigma \sim \varepsilon^{-1}|\ln \varepsilon|^{1-1/\alpha}$ for $\alpha > 1$.

**1. Introduction.** In this paper a mathematically rigorous study of the first exit problem for jump-diffusions driven by small scale Lévy processes with light big jumps is given. The problem under consideration can be outlined as follows. Consider a deterministic dynamical system given by a differential equation $\dot{Y}_t = -U'(Y_t)$ which has a unique asymptotically stable state at the origin 0. We assume that the interval $(-1,1)$ belongs to the domain of attraction of 0, so that the solution trajectories of the deterministic part cannot leave this interval.

The situation becomes different if the dynamical system is perturbed by some noise of small intensity. The stable state becomes meta-stable, and exits from the interval become possible. The probabilistic characteristics of the first exit time, such as its law or the mean value, are determined by the nature of the noise and geometry of the potential $U$.

Received November 2007; revised April 2008.

[1]Supported by the Deutsche Forschungsgemeinschaft (DFG).

Unquestionably, perturbations by Brownian motion are by far the best understood. The first exit problem for small Brownian perturbations was treated in a pioneering work by Kramers [20]. The main mathematical reference on this subject is the book [14] by Freidlin and Wentzell, in which the theory of large deviations for dynamical systems with small Brownian perturbations was developed as the main tool for exit problems. Among many other papers on this subject, we mention a paper by Williams [31], a book [28] by Schuss, and a series of papers by Day [7, 8, 9] and Bovier et al. [3, 4].

In our particular case, the results obtained in very general geometric settings in the references cited above offer a simple explanation. It turns out that the length of the mean exit time is asymptotically given by $e^{\zeta/\varepsilon^2}$, and its logarithmic rate $\zeta$ is determined by the lowest potential barrier a Gaussian particle has to overcome in order to exit. For instance, if $U(-1) < U(1)$, the exit occurs with an overwhelming probability at $-1$, and $\zeta = 2(U(-1) - U(0))$. Moreover, the normalized exit time has a standard exponential law in the limit of small noise $\varepsilon \to 0$.

It is interesting to note that asymptotics of the Gaussian type are also obtained in a situation in which a random Markov perturbation possesses locally infinitely divisible laws with exponential moments of any order, while jump intensity increases and jump size decreases simultaneously with an appropriate rate along with the noise parameter $\varepsilon$ tending to 0. A typical example of such a perturbation is given by a compensated Poisson process with jump size $\varepsilon$, and jump intensity $1/\varepsilon$, see [14], Chapter 5.

The situation becomes quite different for a dynamical system perturbed by heavy tailed Lévy noise. There the big jumps begin to play the major role in the exit time dynamics. If the jump measure of the driving Lévy process has power tails, the mean exit time turns out to behave like a power function of the small noise amplitude. Moreover, the leading term of the average first exit time does not depend on the vertical parameters of the potential's geometry, the heights of potential barriers, but rather on horizontal ones such as the distances between the stable point and the domain boundary. Due to the presence of big jumps, the trajectories of the perturbed dynamical system leave the interval in one big jump, and do not climb up the potential barrier as in the Gaussian case.

Rough large deviation estimates and the asymptotics of the mean first exit time from a domain for a more general class of Markov processes with heavy tailed jumps were first obtained by Godovanchuk [15], whereas a general large deviations theory for Markov processes can be found in the book [29] by Wentzell. In [26], Samorodnitsky and Grigoriu studied the tails of jump-diffusions driven by Lévy processes with heavy (regularly varying) tails. A lot of information about jump-diffusions and models with heavy tails can be found in a book [13] by Embrechts, Klüppelberg and Mikosch and the references therein.

Imkeller and Pavlyukevich [16, 17] recently described the fine asymptotics for the law and moments of exit times for jump diffusions driven by $\alpha$-stable Lévy processes and, more generally, by Lévy processes the jump measure of which has regularly varying tails. These asymptotic properties were used to show metastability properties for Lévy-driven jump diffusions in multi-well potentials (see [16, 18]). The techniques were further generalized to study simulated annealing with time nonhomogeneous jump processes; see [23, 24].

Our recent interest in small noise jump diffusions arose from the acquaintance with the paper [11] by Ditlevsen. In an attempt to model paleoclimatic time series from the Greenland ice core by dynamical systems perturbed by noise, the author discovers an $\alpha$-stable noise signal with $\alpha \approx 1.75$. In his setting, big jumps of the $\alpha$-stable Lévy process are responsible for rapid catastrophic climate changes (the so-called Dansgaard–Oeschger events) observed in the Earth's northern hemisphere during the last glaciation. The appearance of a stable Lévy noise signal can be explained if the observed time series is interpreted as a mesoscopic limit of some more complicated climate dynamics.

Lévy-driven stochastic dynamics has become a popular research field in physics, where stable non-Gaussian Lévy processes are often named *Lévy flights*. We draw the reader's attention to the topical review by Metzler and Klafter [22], where Lévy flights are discussed in detail from a physicists point of view.

The first exit problem (also called Kramers' or barrier crossing problem) is of central importance in the physical sciences. Small noise barrier crossing problems were studied on a physical level of rigor by Ditlevsen [10], Chechkin at al. [5, 6] and Dybiec, Gudowska-Nowak and Hänggi [12]. In particular, Chechkin at al. [5] present numerical experiments that strongly support the theoretical findings of [16, 18].

The relationship between Lévy and Gaussian exit time dynamics circumscribes another problem that has received a considerable deal of attendance in physics applications. The problem was first considered by Mantegna and Stanley in [21] and Koponen [19]. In order to see Gaussian type asymptotic behavior emerge in dynamical systems perturbed by Lévy processes, the authors suggest to either eliminate big jumps, or to make their appearance rare by modifying the jump measure to have exponentially light tails.

In this paper we will study exit times of solutions of the stochastic differential equation

$$(1.1) \qquad dX_t^\varepsilon = -U'(X_{t-}^\varepsilon)\,dt + \varepsilon\,dL_t$$

driven by a Lévy process $L$ of (small) intensity $\varepsilon$ whose jump measure $\nu$ is symmetric, and which contains a nontrivial Gaussian component. The argument and results obtained in [16, 17] for the heavy-tail jump measures when $\nu([u,\infty)) \sim u^{-r}$ with some $r > 0$ show that the exit occurs due to a single big jump of the order $1/\varepsilon$ and, thus, a mean exit time is of the order $\nu(\{|y| \geq 1/\varepsilon\})^{-1} \sim \varepsilon^{-r}$, $\varepsilon \to 0$. Moreover, one can see that the argument of [16] would hold also for Lévy measures with sub-exponential tails $\nu([u,+\infty)) \sim e^{-u^\alpha}$ for very small values $\alpha \ll 1$ leading to mean exit times of the order $e^{1/\varepsilon^\alpha}$. Recalling that Gaussian exits occur on time scales of the order $e^{1/\varepsilon^2}$, one may ask the following: a further reduction of the tail weight can lead to Gaussian dynamics as $\alpha \uparrow 2$, and will Gaussian exit dynamics dominate after crossing the critical index, that is, for $\alpha > 2$? This is the motivating question of this paper.

The result is very surprising. We show that big jumps *always* dominate, independently of how light they are. Looking at the results in more detail, the behavior of exit times encounters a phase transition at the critical value $\alpha = 1$ which marks the transition from *sub-exponential* to *super-exponential* dynamics.

Our arguments leading to the discovery of this transition can be outlined as follows. As in the case for power type tails in [17], for $\varepsilon > 0$ the process $L$ is decomposed into a compound Poisson pure jump part $\eta^\varepsilon$ with jumps of height larger than some critical level $g_\varepsilon$, and a remainder $\xi^\varepsilon$ with jumps not exceeding this bound. The critical threshold $g_\varepsilon$ has to be chosen individually according to whether the jump law possesses sub- or super-exponential tails. In the crucial estimate, one has to show that the exit from the interval $(-1,1)$ around the stable fixed point $0$ before some finite time $T$ while never returning to a small neighborhood of $0$ can occur in two ways. First, the increments of the small jump component $\xi^\varepsilon$ have to exceed certain bounds, and the probability of this scenario can be made small enough by suitable choice of $g_\varepsilon$. Second, the sum of large jumps occurring before time $T$ exceeds the bound $1$. Generally, the analysis reveals that large jumps are responsible for exits irrespective of whether we are in the sub- or the super-exponential regime.

To see the phase transition at $\alpha = 1$, consider the big jumps $W_i$, $i \geq 1$, of the process $L$. The jumps are independent and have the law $\beta_\varepsilon^{-1} \nu|_{[-g_\varepsilon, g_\varepsilon]^c}(\cdot)$ with $\beta_\varepsilon = \nu([-g_\varepsilon, g_\varepsilon]^c) \approx 2\exp(-g_\varepsilon^\alpha)$ being the inverse mean time between big jumps.

The essential part of the proof consists in an asymptotic estimate of the tail probability

$$(1.2) \qquad \mathbf{P}\left(\sum_{i=1}^{k} |\varepsilon W_i| \geq 1\right),$$

which is, roughly speaking, the probability of the exit in no more than $k$ big jumps of $L$, the number $k := k_\varepsilon$ being chosen appropriately. Then the following estimate contributes crucially to the phase transition:

$$(1.3) \quad \mathbf{P}\left(\sum_{i=1}^{k} |\varepsilon W_i| \geq 1\right) \leq c_\varepsilon^k \exp\left(-\inf\left\{\sum_{i=1}^{k} x_i^\alpha : \sum_{i=1}^{k} x_i = \frac{1}{\varepsilon}, x_i \geq 0\right\}\right),$$
$$c_\varepsilon > 0.$$

The phase transition emerges when solving the minimization problem in the exponent of this estimate. Thus, in the case of *sub-exponential tails*, the infimum in the exponent is attained *on the boundary* of the simplex $\{(x_1, \ldots, x_n) : x_i \geq 0, \sum_{i=1}^{n} x_i = \varepsilon^{-1}\}$, namely, at points $x_i = \varepsilon^{-1}$, $x_j = 0$, $j \neq i$.

On the contrary, for *super-exponential tails*, the infimum is attained in the inner point of the simplex, namely, at $x_i = (\varepsilon k)^{-1}$, $1 \leq i \leq k$.

One can say that from the point of view of the optimization technique, the phase transition is due to the switching from *concavity* to *convexity* of the function

$$(1.4) \qquad x \mapsto x^\alpha, \qquad x \geq 0$$

as $\alpha$ increases through 1. So, the surprising behavior of our jump diffusion with exponentially light jumps of degree $\alpha$ can roughly be summarized by this statement: big jumps of the Lévy process govern the asymptotic behavior in the sub- ($\alpha < 1$) as well as in the super- ($\alpha > 1$) exponential regimes, but in the latter one the cumulative action of several large jumps, reminding the climbing of the potential well in the Gaussian regime, becomes important, while for $\alpha < 1$ the biggest jump alone governs the exit.

The material is organized as follows. In Section 2 we explain the setup of the problem and state our main results about the asymptotics of exit times. Section 3 contains our general strategy of estimating the tails of the law of exit times from the knowledge of the tails of the jump measure. The concept of $\varepsilon$-dependent separation of small and large jump parts which already proved to be successful in [17] takes a central role, and is basic for the proof that, apart from a Gaussian part, small jumps do not alter the behavior of solution curves of the unperturbed dynamical system by much. This leaves the role of triggering exits to the large jump part, the contribution of which receives a careful estimation. In the technical Sections 4 and 5 upper and lower bounds for the tails of the exit time laws are derived.

**2. Object of study and main result.** Let $(\Omega, \mathcal{F}, (\mathcal{F}_t)_{t \geq 0}, \mathbf{P})$ be a filtered probability space. We assume that the filtration satisfies the usual hypotheses in the sense of [25], that is, it is right-continuous, and $\mathcal{F}_0$ contains all the $\mathbf{P}$-null sets of $\mathcal{F}$.

For $\varepsilon > 0$ we consider solutions $X^\varepsilon = (X^\varepsilon_t)_{t \geq 0}$ of the one-dimensional stochastic differential equation

$$(2.1) \qquad X^\varepsilon_t(x) = x - \int_0^t U'(X^\varepsilon_{s-}(x))\,ds + \varepsilon L_t, \qquad t \geq 0,\ x \in \mathbb{R},$$

where $L$ is a Lévy process and $U$ is a potential function satisfying assumptions specified in the following. The principal goal of our investigation is the small noise behavior of $X^\varepsilon$, that is, the behavior of the process as $\varepsilon \to 0$. More specifically, our interest is focused on the exit of $X^\varepsilon$ from a neighborhood of the stable attracting point $0$ of the real valued potential function $U$ defined on $\mathbb{R}$. For this reason, besides assuming that $U$ be continuously differentiable with derivative $U'$, we only have to fix some minimal conditions on $U$ concerning its properties in a neighborhood of $0$. We shall work under the following assumption:

(U) $U'(x) = 0$ if $x = 0$, $U'$ is Lipschitz continuous, and $U'(x)x > 0$, $x \in (-1, 1) \setminus \{0\}$.

In particular, the drift $U'$ may vanish at the ends of the interval, $U'(\pm 1) = 0$.

In order to state the conditions our Lévy process $L$ is supposed to satisfy, let us recall that a positive Lebesgue measurable function $l$ is slowly varying at infinity if $\lim_{u \to +\infty} l(\lambda u)/l(u) = 1$ for any $\lambda > 0$. For example, positive constants, logarithms and iterated logarithms are slowly varying functions:

(L1) $L$ has a generating triplet $(d, \nu, \mu)$ with a Gaussian variance $d \geq 0$, an arbitrary drift $\mu \in \mathbb{R}$ and a *symmetric* Lévy measure $\nu$ satisfying the usual condition $\int_{\mathbb{R} \setminus \{0\}} (y^2 \wedge 1)\nu(dy) < \infty$.

(L2) Let $f(u) := -\ln \nu([u, +\infty))$, $u > 0$. Then there is $\alpha > 0$ such that

$$(2.2) \qquad f(u) = u^\alpha l(u), \qquad u \geq 1$$

for some function $l$ slowly varying at $+\infty$.

We say that the Lévy measure $\nu$ has *sub-exponential* or *super-exponential tails* with index $\alpha$ if $0 < \alpha < 1$ or $\alpha > 1$ in (L2), respectively. Typical examples of Lévy processes under consideration are given by symmetric Lévy measures $\nu$ with tails

$$(2.3) \qquad \nu([u, \infty)) = \exp(-u^\alpha), \qquad \alpha > 0,\ u \geq 1.$$

The family of strongly tempered stable processes with the jump measures

$$(2.4) \qquad \nu(dy) = e^{-\lambda |y|^\alpha} |y|^{-1-\beta} \mathbb{I}\{y \neq 0\}\,dy,$$

$$\lambda > 0,\ \beta \in (0, 2),\ \alpha > 0,\ \alpha \neq 1,$$

provides another example.

Since Lévy processes are semimartingales, and due to (U), the stochastic differential equation (2.1) possesses a strong solution defined for all $t \geq 0$. See [1, 25] for the general theory of stochastic integration and [27] for more information on Lévy processes. Moreover, the underlying deterministic equation ($\varepsilon = 0$) given by

$$(2.5) \qquad Y_t(x) = x - \int_0^t U'(Y_s(x))\,ds$$

has a unique solution for any initial value $x \in \mathbb{R}$ and all $t \geq 0$. The position $x = 0$ of the minimum of $U$ is a stable attractor for the dynamical system $Y$, that is, for any $x \in (-1,1)$ we have $Y_t(x) \to 0$ as $t \to \infty$. It is clear that the deterministic solution $Y(x)$ does not leave the interval $[-1,1]$ for initial values $x \in (-1,1)$. The main object of study of this paper is the asymptotic law and the mean value of the first exit time of the jump-diffusion $X^\varepsilon$:

$$(2.6) \qquad \sigma_x(\varepsilon) = \inf\{t \geq 0 : |X_t^\varepsilon(x)| \geq 1\}, \qquad x \in \mathbb{R}.$$

Our main findings are stated in the following theorems.

THEOREM 2.1 (Sub-exponential tails). *Let the jump measure $\nu$ of $L$ be sub-exponential with index $0 < \alpha < 1$. Then for any $\delta > 0$ there is $\varepsilon_0 > 0$ such that for all $0 < \varepsilon \leq \varepsilon_0$ the following inequalities hold uniformly for $t \geq 0$:*

$$(2.7) \quad (1-\delta)\exp(-C_\varepsilon^{1-\delta} t) \leq \inf_{|x| \leq 1-\delta} \mathbf{P}(\sigma_x(\varepsilon) > t) \leq \sup_{|x| \leq 1} \mathbf{P}(\sigma_x(\varepsilon) > t)$$
$$\leq \exp(-\tfrac{1}{2} C_\varepsilon t)$$

*with $C_\varepsilon := \nu((-\frac{1}{\varepsilon}, \frac{1}{\varepsilon})^c) = 2\exp(-f(\frac{1}{\varepsilon}))$. Consequently, for any $|x| < 1$ we have*

$$(2.8) \qquad \lim_{\varepsilon \to 0} f\left(\frac{1}{\varepsilon}\right)^{-1} \ln \mathbf{E}\sigma_x(\varepsilon) = 1.$$

REMARK 2.1. It will be seen from the proof (Section 4.1) that the upper bound in (2.7) holds not only for small $\varepsilon$ but for all $\varepsilon > 0$ and for all symmetric jump measures $\nu$. Moreover, the factor $1/2$ in the exponent of the upper bound (2.7) can be omitted if $C_\varepsilon$ satisfies $C_\varepsilon = \inf_{y \in \mathbb{R}} \nu((\frac{-1-y}{\varepsilon}, \frac{1-y}{\varepsilon})^c)$, which, for instance, holds for unimodal jump measures $\nu$ with the mode 0.

THEOREM 2.2 (Super-exponential tails). *Let the jump measure $\nu$ of $L$ be super-exponential with index $\alpha > 1$. Let $q_\varepsilon$ denote its $\varepsilon$-quantile, $q_\varepsilon := \sup\{u > 0 : \nu([u, \infty)) \geq \varepsilon\}$. Then for any $\delta > 0$ there is $\varepsilon_0 > 0$ such that for all $0 < \varepsilon \leq \varepsilon_0$ the following inequalities hold uniformly for $t \geq 0$:*

$$(2.9) \quad (1-\delta)\exp(-D_\varepsilon^{1-\delta} t) \leq \inf_{|x| \leq 1-\delta} \mathbf{P}(\sigma_x(\varepsilon) > t) \leq \sup_{|x| \leq 1} \mathbf{P}(\sigma_x(\varepsilon) > t)$$
$$\leq (1+\delta)\exp(-D_\varepsilon^{1+\delta} t),$$

where $D_\varepsilon = \exp(-d_\alpha \frac{|\ln \varepsilon|}{\varepsilon q_\varepsilon})$ and $d_\alpha = \alpha(\alpha - 1)^{1/\alpha - 1}$. Consequently, for any $|x| < 1$ we have

$$(2.10) \qquad d_\alpha^{-1} \lim_{\varepsilon \to 0} \frac{\varepsilon q_\varepsilon}{|\ln \varepsilon|} \ln \mathbf{E}\sigma_x(\varepsilon) = 1.$$

It is instructive to compare qualitatively the results of Theorems 2.1 and 2.2 with known results for exit times in the case in which $L$ is a pure Brownian motion, or contains a symmetric jump component with regularly varying tails. We therefore briefly consider the mean exit times of Lévy-driven diffusions of four types. Then the following limiting relations hold and are uniform over all initial points $x$ belonging to a compact subset $K \subset (-1, 1)$:

1. Power tails. Let $L$ be such that $\nu([u, \infty)) = u^{-r}$, $u \geq 1$ for some $r > 0$. Then as was shown in [17], the mean exit time satisfies

$$(2.11) \qquad 2 \lim_{\varepsilon \to 0} \varepsilon^r \mathbf{E}\sigma_x(\varepsilon) = 1.$$

2. Sub-exponential tails. Assume that $L$ is such that for some $\alpha \in (0, 1)$ we have $\nu([u, \infty)) = \exp(-u^\alpha)$, $u \geq 1$. Then Theorem 2.1 easily implies that

$$(2.12) \qquad \lim_{\varepsilon \to 0} \varepsilon^\alpha \ln \mathbf{E}\sigma_x(\varepsilon) = 1.$$

3. Super-exponential tails. Assume that $L$ is such that $\nu([u, \infty)) = \exp(-u^\alpha)$, $u \geq 1$, for some $\alpha \in (1, \infty)$. Then the $\varepsilon$-quantile $q_\varepsilon = |\ln \varepsilon|^{1/\alpha}$ and Theorem 2.2 entails that

$$(2.13) \qquad \frac{1}{\alpha}(\alpha - 1)^{1 - 1/\alpha} \lim_{\varepsilon \to 0} \varepsilon |\ln \varepsilon|^{1/\alpha - 1} \ln \mathbf{E}\sigma_x(\varepsilon) = 1.$$

4. Gaussian diffusion. Assume that $L$ possesses the characteristic triplet $(1, 0, 0)$, that is, $L$ is a standard one-dimensional Brownian motion. Then the mean exit time depends on the height of the potential barrier at the interval ends, and

$$(2.14) \qquad \tfrac{1}{2}(U(-1) \wedge U(1))^{-1} \lim_{\varepsilon \to 0} \varepsilon^2 \ln \mathbf{E}\sigma_x(\varepsilon) = 1.$$

First we note that cases 1 and 2 mathematically do not differ by much, since the mean exit time can be expressed by the same formula $\mathbf{E}\sigma_x(\varepsilon) \sim (2\nu([\frac{1}{\varepsilon}, \infty)))^{-1}$. The gaps between 2 and 3, and 3 and 4 are much more surprising. The logarithmic rate of the expected first exit time drastically changes its behavior in the super-exponential case: jump lightness influences the mean exit time in a rather insignificant way. Even more surprising is the fact that we do not obtain Gaussian asymptotics even for light tails with $\alpha \geq 2$. This is underlined in an intriguing way through the form of the pre-factors: in the cases of perturbations with jumps they only depend on the distance between the stable equilibrium 0 and the interval boundary, whereas in the Gaussian case the heights of the potential barriers come into play.

To say more, the phase transition between 3 and 4 shows that the transition to Gaussian dynamics is impossible with Lévy perturbations of the type $\varepsilon L$, that is, by scaling only *sizes* of jumps and not their *intensity*. However, as we already mentioned in the Introduction, Gaussian exit times can be obtained if we couple the *size* and *intensity* of jumps.

Finally, we apply the tools developed for Theorems 2.2 and 2.1 to study another class of perturbations with *bounded jumps*, which leads to exit times of the order $\nu([\frac{1}{\varepsilon},\infty))^{-a}$ for arbitrary $a > 0$, $\nu$ being a symmetric sub-exponential Lévy measure with $\alpha \in (0,1)$.

For any $\theta > 0$, consider a Lévy process with *bounded jumps* $L^{\varepsilon,\theta}$ with a characteristic triplet $(d, \nu^{\varepsilon,\theta}, \mu)$, where

$$\nu^{\varepsilon,\theta} = \nu|_{[-\theta/\varepsilon, \theta/\varepsilon]}, \qquad \theta > 0, \tag{2.15}$$

$d \geq 0$ and $\mu \in \mathbb{R}$. The corresponding jump-diffusion $X^{\varepsilon,\theta}$ is a strong solution of (2.1) with $L^{\varepsilon,\theta}$ instead of $L$. In this setting, the jumps of the process $\varepsilon L^{\varepsilon,\theta}$ and, hence, of $X^{\varepsilon,\theta}$ are bounded by $\varepsilon$-independent value $\theta$, which makes impossible the exit of $X^{\varepsilon,\theta}$ from a neighborhood of 0 in a single big jump if $\theta < 1$.

THEOREM 2.3 (Bounded sub-exponential tails). *For $\alpha \in (0,1)$ and $\theta > 0$, let $L^{\varepsilon,\theta}$ have the jump measure $\nu^{\varepsilon,\theta}$. Then for any $\delta > 0$ there is $\varepsilon_0 > 0$ such that for all $0 < \varepsilon \leq \varepsilon_0$ the following inequalities hold uniformly for $t \geq 0$:*

$$(1-\delta)\exp(-C_{\varepsilon,\theta}^{1-\delta} t) \leq \inf_{|x| \leq 1-\delta} \mathbf{P}(\sigma_x(\varepsilon) > t) \leq \sup_{|x| \leq 1} \mathbf{P}(\sigma_x(\varepsilon) > t)$$
$$\leq (1+\delta)\exp(-C_{\varepsilon,\theta}^{1+\delta} t) \tag{2.16}$$

*with $C_{\varepsilon,\theta} := \nu([\frac{1}{\varepsilon},\infty))^{\phi(\theta)}$ and*

$$\phi(\theta) := \left[\frac{1}{\theta}\right]\theta^\alpha + \left(1 - \left[\frac{1}{\theta}\right]\theta\right)^\alpha$$
$$= \inf\left\{\sum_{i=1}^k x_i^\alpha : \sum_{i=1}^k x_i = 1, x_i \in [0,\theta], k \geq 1/\theta\right\}. \tag{2.17}$$

*Consequently, for any $|x| < 1$ we have*

$$\lim_{\varepsilon \to 0}\left(\phi(\theta) f\left(\frac{1}{\varepsilon}\right)\right)^{-1} \ln \mathbf{E}\sigma_x(\varepsilon) = 1. \tag{2.18}$$

In particular, if $\nu([u,\infty)) = \exp(-u^\alpha)$, $u \geq 1$ with $\alpha \in (0,1)$ and $\theta = \frac{1}{k}$, $k \geq 1$, we have $\phi(\theta) = k^{1-\alpha}$, and Theorem 2.3 yields

$$\lim_{\varepsilon \to 0} \varepsilon^\alpha \ln \mathbf{E}\sigma_x(\varepsilon) = k^{1-\alpha}. \tag{2.19}$$

Finally, we note that an interested reader can find more details on the cases $\alpha = 0$ (in particular, slowly varying $f(u) = -\ln \nu([u, \infty))$) and $\alpha = +\infty$ (in particular, jump measures with bounded support), as well as on the critical exponential case $\alpha = 1$ in [30].

## 3. Main tools for the proof.

3.1. *Key elements of the proof.* First, the estimates of Theorems 2.1, 2.2 and 2.3 will follow from essentially elementary but very general inequalities which allow to determine the bounds for the probability distribution function of the first exit time in terms of the process's dynamics on the fixed time intervals (Lemma 3.1). Whereas the upper estimate of the probability $\mathbf{P}(\sigma_x > t)$, $t \geq 0$, in terms of $\mathbf{P}(\sigma_x \leq T)$, $T > 0$ fixed, is well known, the lower estimate requires consideration of the event $\{\sigma_x^* \leq T\}$, $\sigma_x^*$ being a *non-Markovian* time of the *last exit* from some domain.

Next, in Section 3.3 we decompose the driving Levy process into small and big jump parts. We show that on the event $\{\sigma_x^* \leq T\}$, if the exit from the interval $(-1, 1)$ occurs after the $k$th big jump, then either the small jump part makes a big deviation on some short time interval, or $k$ big jumps make up a sequence with short interjump times.

The exponential bound for the probability of big deviation of the small jump process is obtained in Section 3.4. Further, we derive an exponential tail estimate for sums of big jumps expressed in terms of a multivariate minimization problem with constraints.

With these tools, in Sections 4 and 5 we carefully choose the critical $\varepsilon$-dependent threshold to separate big and small jumps, determine the $\varepsilon$-dependent number of jumps which contribute mainly to the exit, and finally prove the main results of this paper.

3.2. *Estimates on short time intervals.* Let $X = (X_t(x))_{t \geq 0}$ be a time homogeneous Markov process starting in $x \in \mathbb{R}$ whose sample paths are right-continuous and have left limits. Consider an interval $I \subset \mathbb{R}$, a subinterval $J \subset I$, and consider the (Markovian) *first exit* time $\sigma_x := \inf\{t \geq 0 : X_t(x) \notin I\}$, and the (non-Markovian) time $\sigma_x^* := \sup\{t < \sigma : X_t(x) \in J\}$ which marks the *start of the exit*.

The following lemmas allow to estimate the law of $\sigma_x$ from the dynamics of the process $X$ on relatively short time intervals. For the sake of simplicity of notation, we sometimes omit the subscript $x$ in expressions containing the times $\sigma_x$ and $\sigma_x^*$.

LEMMA 3.1. *Let $C$ and $T$ be positive real numbers such that $CT < 1$.*

1. If $\inf_{x\in I} \mathbf{P}(\sigma_x \leq T) \geq CT$, the following estimate from above for $t \geq 0$ holds:

$$\sup_{x\in I} \mathbf{P}(\sigma_x > t) \leq (1-CT)^{-1} \exp(-Ct). \tag{3.1}$$

2. If $\sup_{x\in J} \mathbf{P}(\sigma_x^* \leq T) \leq CT$, the following estimate from below for $t \geq 0$ holds:

$$\inf_{x\in J} \mathbf{P}(\sigma_x > t) \geq (1-CT) \exp\left(\frac{\ln(1-CT)}{T} t\right). \tag{3.2}$$

PROOF. The proof of part 1 is a straightforward application of the strong Markov property and time homogeneity of $X$. In fact, choose an arbitrary $t > 0$ and let $k := [\frac{t}{T}]$. Then for any $x \in I$ we obtain the following chain of inequalities:

$$\mathbf{P}(\sigma_x > t) \leq \mathbf{P}(\sigma_x > kT) \leq \left(\sup_{x\in I} \mathbf{P}(\sigma_x > T)\right)^k$$
$$\leq (1-CT)^{t/T-1} \leq (1-CT)^{-1} \exp(-Ct). \tag{3.3}$$

In order to use similar arguments to prove part 2, we need to define a sequence of stopping times. Let $T_J^0 := 0$ and for any $n \geq 1$ let $T_J^n := \inf\{t: t \geq T_J^{n-1} + T, X_t(x) \in J\}$. Obviously $\{\sigma_x \leq T_J^1\} = \{\sigma_x^* \leq T\}$ holds for any $x \in J$ and, moreover, $T_J^n \geq nT$ for any $n \geq 1$. Again choose an arbitrary $t > 0$ and let $k := [\frac{t}{T}]$. Then for any $x \in J$ we have

$$\mathbf{P}(\sigma_x > t) \geq \mathbf{P}(\sigma_x > T_J^{k+1}) \geq \left(\inf_{x\in J} \mathbf{P}(\sigma_x > T_J^1)\right)^{k+1}$$
$$\geq (1-CT)^{t/T+1} \geq (1-CT) \exp\left(\frac{\ln(1-CT)}{T} t\right). \quad \square \tag{3.4}$$

3.3. *Decomposition into small and large jump parts.* In our separation of the jump part of the Lévy process $L$ into a component for small and one for large jumps the latter will turn out to be a compound Poisson process. This makes large jumps relatively easily amenable to an individual investigation. Suppose that $g > 0$ is a *cutoff* height. We shall leave a particular choice of $g$ to later parts of this paper, and for the moment use the cutoff height to define the $g$-dependent decomposition

$$L = \xi + \eta \tag{3.5}$$

with jump measures $\nu_\xi = \nu|_{[-g,g]}$ and $\nu_\eta = \nu|_{[-g,g]^c}$. The resulting independent Lévy processes $\eta$ and $\xi$ possess generating triplets $(0, \nu_\eta, 0)$ and $(d, \nu_\xi, \mu)$. For the compound Poisson part $\eta$ and $k \geq 1$ we denote by $S_k$ the arrival times of jumps ($S_0 = 0$), by $\tau_k = S_k - S_{k-1}$ its interjump periods, and by $W_k$ the respective jump sizes, and note that $\beta = \nu_\eta(\mathbb{R})$ expresses its jump frequency, that is, the inverse expected interjump time. Finally, we denote by $N_t = \sup\{k \geq 0 : S_k \leq t\}$ the number of jumps until time $t, t \geq 0$.

Lemma 3.1 reduces the main task of the proof of Theorems 2.1 and 2.2 to finding an appropriate $T > 0$ and estimating the probabilities to exit before $T$ and to start an exit from a subinterval before $T$. For technical reasons, we have to reduce the interval $I = (-1, 1)$ a bit, and study exits from this subinterval. So for some $0 < \delta < \frac{1}{2}$ let $I_\delta^- := (-1 + \delta, 1 - \delta)$, and $\sigma_x^- := \inf\{t : X_t^\varepsilon(x) \notin I_\delta^-\}$. Now take $J_\delta = [-\delta, \delta]$ as a subinterval of $I_\delta^-$ and, according to Lemma 3.1, consider $\sigma_x^* := \sup\{t \leq \sigma_x^- : X_t^\varepsilon(x) \in J_\delta\}$. The following auxiliary estimates intend to control the probability of $\{\sigma_x^* < T\}$ through finding a covering by sets of sufficiently small probability.

LEMMA 3.2.  (i) Let $0 < \delta < \frac{1}{2}$ and $\varepsilon > 0$ be such that $\varepsilon g < \delta$, and $x \in I_\delta^-$. Let $m := \inf_{y \in I_\delta^- \setminus J_\delta} |U'(y)|$ and $\hat{T} > \frac{1}{m}$. Then for any $T > 0$ and $t > T + \hat{T}$, the following estimate holds:

$$\{\sigma_x^* < T, \sigma_x^- \geq t\} \subseteq \{t - S_{N_t} < \hat{T}\} \cup \left\{\sup_{r \leq \hat{T}} \varepsilon |\xi_{t-\hat{T}+r} - \xi_{t-\hat{T}}| \geq m\hat{T} - 1\right\}.$$
(3.6)

(ii) For any $x \in I$ and $T \in [0, \sigma_x]$, we have the following estimate:

$$|X_T^\varepsilon(x)| \leq |x| + \left(\sup_{t \leq T} - \inf\right) \varepsilon L_t.$$
(3.7)

PROOF. (i) Choose $T > 0$, $\hat{T} > \frac{1}{m}$ and $t \geq T + \hat{T}$ arbitrarily. Consider the event $A := \{\sigma_x^* < T, \sigma_x^- \geq t, t - S_{N_t} \geq \hat{T}\}$. It is sufficient to show that on $A$ we have

$$\sup_{r \leq \hat{T}} \varepsilon |\xi_{t-\hat{T}+r} - \xi_{t-\hat{T}}| \geq m\hat{T} - 1.$$
(3.8)

First of all, by definition, on $A$ the process $\varepsilon L$ cannot make jumps larger than $\varepsilon g$ during the time period $[t - \hat{T}, t]$ and does not leave $I_\delta^- \setminus J_\delta$ during the time period $[T, t]$. Further, by choice of $T$ and $\hat{T}$, we have $[t - \hat{T}, \hat{T}] \subseteq [T, t]$ and $\varepsilon g \leq \delta$. Thus, $X^\varepsilon(x)$ does not change its sign during the time period $[t - \hat{T}, \hat{T}]$. By definition, we have $A \subseteq \{X_t^\varepsilon(x) \neq 0\}$. Hence, it is sufficient to consider the cases $A \cap \{X_t^\varepsilon(x) > 0\}$ and $A \cap \{X_t^\varepsilon(x) < 0\}$ separately. On the former we have

$$0 < X_t^\varepsilon(x) = X_{t-\hat{T}}^\varepsilon(x) - \int_{t-\hat{T}}^t U'(X_{s-}^\varepsilon(x))\,ds + \varepsilon(L_t - L_{t-\hat{T}})$$
(3.9)
$$\leq 1 - m\hat{T} + \sup_{s \leq \hat{T}} \varepsilon |\xi_{t-\hat{T}+s} - \xi_{t-\hat{T}}|.$$

Analogously, on $A \cap \{X_t^\varepsilon(x) < 0\}$ we may estimate

$$0 > -1 + m\hat{T} - \sup_{s \leq \hat{T}} \varepsilon |\xi_{t-\hat{T}+s} - \xi_{t-\hat{T}}|. \tag{3.10}$$

This completes the proof of part (i).

(ii) For $x \in I$, let $\varrho_x := \inf\{t \in [0,T] : X_s^\varepsilon(x) > 0 \text{ for all } s \in [t,T]\}$. By construction, $X_t^\varepsilon(x) \geq 0$ and $U'(X_t^\varepsilon(x)) \geq 0$ for any $t \in (\varrho_x, T)$. Thus,

$$\begin{aligned}
X_T^\varepsilon(x) &= X_{\varrho_x}^\varepsilon(x) - \int_{\varrho_x}^T U'(X_{t-}^\varepsilon(x))\, dt + \varepsilon(L_T - L_{\varrho_x}) \\
&\leq X_{\varrho_x}^\varepsilon(x) + \varepsilon(L_T - L_{\varrho_x}) \\
&\leq \begin{cases} x + \sup_{t \leq T} \varepsilon L_t, & \text{if } \varrho_x = 0 \\ \varepsilon(L_T - L_{\varrho_x-}), & \text{if } \varrho_x > 0 \end{cases} \leq |x| + \left(\sup - \inf_{t \leq T}\right)\varepsilon L_t.
\end{aligned} \tag{3.11}$$

Analogously, we have $X_T^\varepsilon(x) \geq -|x| - (\sup - \inf_{t \leq T})\varepsilon L_t$. □

COROLLARY 3.1. (i) Let $0 < \delta < \frac{1}{2}$ and $\varepsilon > 0$ be such that $\varepsilon g < \delta$, and $x \in I_\delta^-$. Let $m := \inf_{y \in I_\delta^- \setminus J_\delta} |U'(y)|$ and $T = \frac{2}{m}$. Then for any $k \geq 1$, the following inclusions hold:

$$\{\sigma_x^* < T, \sigma_x^- \geq S_k\} \subseteq \{\sigma_x^* < T, \sigma_x^- \geq 2Tk \wedge S_k\} \subseteq \chi_k \cup \bigcap_{i=1}^k \{\tau_i \leq 2T\} \tag{3.12}$$

with

$$\chi_k = \bigcup_{i=0}^{k-1} \left\{\sup_{t \leq T} \varepsilon |\xi_{S_i+T+t} - \xi_{S_i+T}| \geq 1\right\}. \tag{3.13}$$

(ii) For any $x \in I$, $k \geq 1$ and $T \in [0, \sigma_x]$, the following estimate holds:

$$\sup_{t < S_k \wedge T} |X_t| \leq |x| + \sum_{i=1}^{k-1} |\varepsilon W_i| + 2\sup_{t \leq T} |\varepsilon \xi_t|. \tag{3.14}$$

PROOF. (i) Obviously, it suffices to prove the second inclusion in (3.12). We set $\hat{T} = T = \frac{2}{m}$, and let $\hat{\chi}_i =: \{\sup_{t \leq T} \varepsilon |\xi_{S_i+T+t} - \xi_{S_i+T}| \geq 1\}$, $i \geq 0$, and $x \in I_\delta^-$. Then we notice that

$$\begin{aligned}
&\{\sigma_x^* < T, \sigma_x^- \geq 2Tk \wedge S_k\} \\
&\subseteq \{\sigma_x^* < T, \sigma_x^- \geq S_k\} \cup \{\sigma_x^* < T, \sigma_x^- \geq 2Tk, S_k > 2Tk\}.
\end{aligned} \tag{3.15}$$

For $k \geq 1$, the event $\{S_k > 2Tk\}$ implies that $\tau_i > 2T$ and $S_{i-1} + 2T \leq 2Tk$ for at least one $i$, $1 \leq i \leq k$, and, therefore, using the equality $\{\tau_i > 2T\} =$

$\{N_{S_{i-1}+2T} = i-1\}$ and applying Lemma 3.2(i) with $t = S_{i-1}+2T$, we obtain

$$\{\sigma_x^* < T, \sigma_x^- \geq 2Tk, S_k > 2Tk\}$$
$$\subseteq \bigcup_{i=1}^{k}(\{\sigma_x^* < T, \sigma_x^- \geq 2Tk\} \cap \{\tau_i > 2T, S_{i-1}+2T \leq 2Tk\})$$
(3.16)
$$\subseteq \bigcup_{i=1}^{k}(\{\sigma_x^* < T, \sigma_x^- \geq S_{i-1}+2T\} \cap \{\tau_i > 2T\})$$
$$\subseteq \bigcup_{i=1}^{k}((\hat{\chi}_{i-1} \cup \{S_{i-1}+2T - S_{N_{S_{i-1}+2T}} < T\}) \cap \{N_{S_{i-1}+2T} = i-1\})$$
$$\subseteq \bigcup_{i=1}^{k}\hat{\chi}_{i-1} \subseteq \chi_k.$$

Next we prove the inclusions

(3.17) $\qquad \{\sigma_x^* < T, \sigma_x^- \geq S_i\} \cap \hat{\chi}_{i-1}^c \subseteq \{\tau_i \leq 2T\}, \qquad 1 \leq i \leq k,$

and taking intersections of their right- and left-hand sides over, $i$ we obtain the proper covering for the set $\{\sigma_x^* < T, \sigma_x^- \geq S_k\}$. Consider the decomposition

(3.18)
$$\{\sigma_x^* < T, \sigma_x^- \geq S_i\} \cap \hat{\chi}_{i-1}^c$$
$$\subseteq (\{\sigma_x^* < T, \sigma_x^- \geq S_{i-1}+2T\} \cap \hat{\chi}_{i-1}^c) \cup \{S_i \leq \sigma_x^- < S_{i-1}+2T\}.$$

The second set in the previous formula is obviously contained in $\{\tau_i \leq 2T\}$, while to study the first one we apply again Lemma 3.2(i) with $t = S_{i-1}+2T$ to obtain

(3.19)
$$\{\sigma_x^* < T, \sigma_x^- \geq S_{i-1}+2T\} \cap \hat{\chi}_{i-1}^c$$
$$\subseteq (\{S_{i-1}+T - S_{N_{S_{i-1}+2T}} < 0\} \cup \hat{\chi}_{i-1}) \cap \hat{\chi}_{i-1}^c$$
$$\subseteq \{N_{S_{i-1}+2T} > i-1\} = \{\tau_i \leq 2T\}.$$

(ii) The second part follows from the estimate

(3.20)
$$\left(\sup_{t<S_k \wedge T} - \inf_{t<S_k \wedge T}\right)\varepsilon L_t \leq \left(\sup_{t<S_k} - \inf_{t<S_k}\right)\varepsilon \eta_t + \left(\sup_{t<T} - \inf_{t<T}\right)\varepsilon \xi_t$$
$$\leq \sum_{i=1}^{k-1}|\varepsilon W_i| + 2\sup_{t \leq T}|\varepsilon \xi_t|. \qquad \square$$

3.4. *Estimates of the small jump process $\xi$.* In this subsection we shall give some estimates for the tails of the maximal fluctuation of the small jump component in the decomposition of noise derived in the previous subsection. In the statement we intend to keep the dependence on the parameters as general as possible, and as explicit as necessary later.

LEMMA 3.3. *Let $\nu \neq 0$ be a symmetric Lévy measure, $b \geq 0$ and $\rho \in \mathbb{R}$. For $g \geq 1$ let $\xi = (\xi_t)_{t \geq 0}$ denote the Lévy process defined by the characteristic triple $(b, \nu|_{[-g,g]}, \rho)$. Then for any $\delta > 0$ there exists $u_0 > 0$ such that, for $T > 0, g \geq 1, f \geq g$ satisfying $\frac{f}{Tg} \geq u_0$, the following estimate holds:*

$$(3.21) \qquad \mathbf{P}\left(\sup_{t \leq T} |\xi_t| > f\right) \leq \exp\left(-(1-\delta)\frac{f}{g}\ln\frac{f}{gT}\right).$$

REMARK 3.1. In particular, if we parameterize $g = g_\varepsilon$, $f = f_\varepsilon$ and $T = T_\varepsilon$ and assume that $\frac{f_\varepsilon}{g_\varepsilon T_\varepsilon} \to \infty$ as $\varepsilon \to 0$, then for every $\delta > 0$ there exists $\varepsilon_0 > 0$, such that (3.21) holds for any $0 < \varepsilon \leq \varepsilon_0$.

PROOF OF LEMMA 3.3. First let us consider the case $\rho = 0$. In this case, for any $g > 0$, the process $\xi$ is a martingale. For any $h > 0$, the reflection principle for symmetric Lévy processes and the Chebyshev inequality applied to the exponent of $\xi_T$ yield

$$(3.22) \qquad \mathbf{P}\left(\sup_{t \leq T} |\xi_t| > f\right) \leq 4\mathbf{P}(\xi_T > f) \leq 4e^{-hf}\mathbf{E}e^{h\xi_T}.$$

The Lévy measure of $\xi$ has bounded support. Thus, the analytic extension of the characteristic function of $\xi$ can be used to estimate $\mathbf{E}e^{h\xi_T}$. Denote $m := \int_{\mathbb{R}}(1 \wedge y^2)\nu(dy) < \infty$, and let $h := \frac{1}{g}\ln\frac{f}{gT}$. Now recalling that $g \geq 1$, we can choose $u_0$ large enough, such that $\frac{b}{2}h^2 \leq \frac{b}{2}(\ln\frac{f}{gT})^2 \leq m\frac{f}{gT}$ for any $\frac{f}{gT} \geq u_0$. The extension of the characteristic function of $\xi$ yields the chain of inequalities

$$(3.23) \quad \begin{aligned} \mathbf{E}e^{h\xi_T} &= \exp\left[\frac{b}{2}h^2 T + T\int_{|y| \leq g}(e^{hy} - 1 - hy\mathbb{I}\{|y| < 1\})\nu(dy)\right] \\ &= \exp\left[\frac{b}{2}h^2 T + T\int_{|y| \leq g}\left(hy\mathbb{I}\{|y| \geq 1\} + \sum_{k=2}^{\infty}\frac{(hy)^k}{k!}\right)\nu(dy)\right] \\ &\leq \exp\left[\frac{b}{2}h^2 T + T\int_{|y| \leq g}\left(hg(1 \wedge y^2) + \sum_{k=2}^{\infty}\frac{(hg)^k}{k!}(1 \wedge y^2)\right)\nu(dy)\right] \\ &\leq \exp\left[\frac{mf}{g} + Tm\exp(hg)\right] = \exp\left[\frac{2mf}{g}\right]. \end{aligned}$$

The statement of the lemma for $\rho = 0$ follows immediately from (3.22) for sufficiently large $u_0$, such that $(2m + \ln 4)/\ln(\frac{f}{gT}) < \delta$ for $\frac{f}{gT} \geq u_0$.

If $\rho \neq 0$, we apply the previous argument to the symmetric martingale $(\xi_t - \rho t)_{t \geq 0}$ and use the estimate

$$\mathbf{P}\left(\sup_{t \leq T} |\xi_t| > f\right) \leq \mathbf{P}\left(\sup_{t \leq T} |\xi_t - \rho t| > f - |\rho|T\right)$$

(3.24)
$$\leq \mathbf{P}\left(\sup_{t\leq T}|\xi_t - \rho t| > (1-\delta')f\right),$$

which holds for any $\delta' > 0$ and sufficiently large $\frac{f}{T} \geq \frac{f}{gT}$. □

3.5. *Tail estimates for the sum of big jumps of $\eta$.* In this crucial subsection we shall give tail estimates for finite sums of jump heights by exponential rates depending on sums of logarithmic tails of the jump laws. The asymptotics of the exit times considered will later be seen to depend on the convexity properties of this sum of logarithmic tails.

Let again $\nu \neq 0$ be a symmetric Lévy measure. For $g \geq 1$, let $\eta = (\eta_t)_{t\geq 0}$ be the Lévy process defined by the characteristic triple $(0, \nu|_{[-g,g]^c}, 0)$. It is clear that $\eta$ is a compound Poisson process. Let $W_k$, $k \geq 1$, denote its jump sizes. The random variables $W_k$ are i.i.d. and satisfy $|W_k| \geq g$.

For $u > 0$ denote $f(u) = -\ln \nu([g \vee u, \infty))$, with the convention $\ln 0 = -\infty$. Let $\beta := \nu([-g,g]^c) = 2\exp(-f(g))$.

LEMMA 3.4. *For $k \geq 1$, let $r$ and $g$ be such that $r > kg$. Then for any $\delta \in (0,1)$ such that $(1-\delta)r > kg$ the following estimate holds:*

$$\mathbf{P}\left(\sum_{i=1}^{k}|W_i| > r\right)$$

(3.25)
$$\leq \frac{2^k}{\beta^k}\left(2 + \frac{\ln r - \ln g}{\ln(1+\delta)}\right)^k$$

$$\times \exp\left(-\inf\left\{\sum_{i=1}^{k}f(x_i): \sum_{i=1}^{k}x_i = (1-\delta)r, x_i \in [g,r]\right\}\right).$$

PROOF. Denote $A := \{(x_1, \ldots, x_k) \in [g,r]^k : \sum_{i=1}^{k} x_i \geq r\}$. We have

(3.26)  $\mathbf{P}\left(\sum_{i=1}^{k}|W_i| > r\right) \leq k\mathbf{P}(|W_1| \geq r) + \mathbf{P}((|W_1|, \ldots, |W_k|) \in A).$

The first summand can be estimated as

$$\mathbf{P}(|W_1| \geq r)$$

(3.27)
$$\leq 2\beta^{-1}\exp(-f(r))$$
$$= 2^k \beta^{-k} \exp(-(f(r) + (k-1)f(g)))$$
$$\leq 2^k \beta^{-k} \exp\left(-\inf\left\{\sum_{i=1}^{k}f(x_i): \sum_{i=1}^{k}x_i \geq (1-\delta)r, x_i \in [g,r]\right\}\right).$$

To estimate the second summand, we cover the set $A$ by a union of parallelepipeds of a special form. Let $M := [\frac{\ln r - \ln g}{\ln(1+\delta)}]$ and consider the set of points

(3.28) $$S := \{s_i, 0 \leq i \leq M\}, \qquad s_i = (1+\delta)^i g$$

and note that $s_M \leq r$ and $(1+\delta)s_M > r$. Consider $(M+1)^k$ parallelepipeds of the type

(3.29) $$P_{t_1,\ldots,t_k} := [t_1, (1+\delta)t_1] \times \cdots \times [t_k, (1+\delta)t_k], \qquad t_1, \ldots, t_k \in S.$$

Obviously, the union of these $(M+1)^k$ parallelepipeds covers the cube $[g, r]^k$, and thus the set $A$. Let $N$ denote the smallest covering of $A$ by these parallelepipeds. If some $P_{t_1,\ldots,t_k} \in N$, that is, $P_{t_1,\ldots,t_k} \cap A \neq \varnothing$, then

(3.30) $$r \leq \max_{(x_1,\ldots,x_k) \in P_{t_1,\ldots,t_k}} \sum_{i=1}^k x_i = (1+\delta) \sum_{i=1}^k t_i$$

and, thus, $\sum_{i=1}^k t_i \geq (1+\delta)^{-1} r \geq (1-\delta)r$. Then

$$\mathbf{P}((|W_1|, \ldots, |W_k|) \in A)$$

$$\leq \sum_{P_{t_1,\ldots,t_k} \in N} \prod_{j=1}^k \mathbf{P}(|W_j| \in [t_j, (1+\delta)t_j])$$

$$\leq (M+1)^k \max_{P_{t_1,\ldots,t_k} \in N} \prod_{j=1}^k \mathbf{P}(|W_j| \in [t_j, (1+\delta)t_j])$$

(3.31) $$\leq (M+1)^k \max_{P_{t_1,\ldots,t_k} \in N} \prod_{j=1}^k \mathbf{P}(|W_1| \geq t_j)$$

$$\leq 2^k (M+1)^k \beta^{-k} \max_{P_{t_1,\ldots,t_k} \in N} \exp\left(-\sum_{j=1}^k f(t_j)\right)$$

$$\leq 2^k (1+M)^k \beta^{-k}$$

$$\times \exp\left(-\inf\left\{\sum_{j=1}^k f(x_j) : \sum_{j=1}^k x_j \geq (1-\delta)r, x_j \in [g, r]\right\}\right).$$

The monotonicity of $f$ and (3.27) lead to

$$\mathbf{P}\left(\sum_{i=1}^k |W_i| > r\right)$$

$$\leq 2^k (k + (1+M)^k) \beta^{-k}$$

$$(3.32) \qquad \times \exp\left(-\inf\left\{\sum_{i=1}^{k} f(x_i) : \sum_{i=1}^{k} x_i \geq (1-\delta)r, x_i \in [g, r]\right\}\right)$$

$$= 2^k(k + (1+M)^k)\beta^{-k}$$

$$\times \exp\left(-\inf\left\{\sum_{i=1}^{k} f(x_i) : \sum_{i=1}^{k} x_i = (1-\delta)r, x_i \in [g, r]\right\}\right).$$

Finally, the elementary inequality $k + (1+M)^k \leq (2+M)^k$, $k \geq 1$, $M \geq 0$, applied to the prefactor completes the estimation. □

3.6. *A simple minimization problem.* Later in Sections 5.2.2 and 5.3.2 we will apply Lemma 3.4 to estimate the tails of sums of big jumps of the process $L$. We shall use the following result. Let $k \geq 1$ and $a > 0$. Then for $0 < \alpha \leq 1$ and $\theta > 0$, we have

$$(3.33) \qquad \inf\left\{\sum_{i=1}^{k} x_i^\alpha : \sum_{i=1}^{k} x_i = a, x_i \in [0, \theta a], k \geq \frac{1}{\theta}\right\}$$
$$= \left[\frac{1}{\theta}\right](\theta a)^\alpha + \left(a - \left[\frac{1}{\theta}\right]\theta a\right)^\alpha$$

and, in particular, for $1 \leq \theta \leq +\infty$,

$$(3.34) \qquad \inf\left\{\sum_{i=1}^{k} x_i^\alpha : \sum_{i=1}^{k} x_i = a, x_i \in [0, \theta a]\right\} = a^\alpha.$$

On the other hand, for $\alpha \geq 1$, we have

$$(3.35) \qquad \inf\left\{\sum_{i=1}^{k} x_i^\alpha : \sum_{i=1}^{k} x_i = a, x_i \geq 0\right\} = k\left(\frac{a}{k}\right)^\alpha.$$

A straightforward application of the method of Lagrangian multipliers implies that the local extremum of the function to minimize is attained for $x_1 = \cdots = x_k = a/k$. Since for $\alpha \geq 1$ this extremum is a local (and global) minimum, the equality (3.35) follows.

In the sub-exponential case, the point determined above is a local maximum and, therefore, the minimum should be looked for on the boundary of the domain. This leads to the equalities (3.34) and (3.33).

The essentially different solutions to the minimization problems above are due to *convex*, respectively, *concave* behavior of the mapping $x \mapsto x^\alpha$ for $\alpha \in (0, 1]$, respectively, $\alpha \geq 1$.

**4. The upper bounds.** In this section we employ the first part of Lemma 3.1 in order to deduce upper bounds for the law of exit times presented in Theorems 2.1, 2.2 and 2.3. This is done in separate arguments in the sub-, super-exponential cases and the case of bounded jumps.

4.1. *Sub-exponential tails. Proof of Theorem 2.1, upper bound.* For the choice $g_\varepsilon := \frac{1}{2\varepsilon}$, $\varepsilon > 0$, let us consider a decomposition $L = \xi^\varepsilon + \eta^\varepsilon$ as in Section 3.3. Let $z^\varepsilon(x)$ be a solution of the corresponding stochastic differential equation driven by small jumps, namely,

$$(4.1) \qquad z_t^\varepsilon(x) = x - \int_0^t U'(z_{s-}^\varepsilon(x))\,ds + \varepsilon \xi_t^\varepsilon.$$

By construction, we have $X_{\tau_1}^\varepsilon(x) = z_{\tau_1}^\varepsilon(x) + \varepsilon W_1$. Thus, for any $\varepsilon > 0$, $T > 0$, and $x \in I$, the following estimate holds:

$$(4.2) \qquad \begin{aligned} \{\sigma_x < T\} &\supseteq \{\tau_1 < T, X_{\tau_1}^\varepsilon(x) \notin I\} \\ &= \{\tau_1 < T, \varepsilon W_1 \notin (-1 - z_{\tau_1}^\varepsilon(x), 1 - z_{\tau_1}^\varepsilon(x))\}. \end{aligned}$$

Since $\tau_1$, $W_1$ and $z^\varepsilon(x)$ are independent, we get

$$(4.3) \qquad \begin{aligned} \mathbf{P}(\sigma_x < T) &\geq \mathbf{P}(\tau_1 < T) \inf_{y \in \mathbb{R}} \mathbf{P}(\varepsilon W_1 \notin (-1-y, 1-y)) \\ &= (1 - e^{-\beta_\varepsilon T})\beta_\varepsilon^{-1} \inf_{y \in \mathbb{R}} \nu_\eta^\varepsilon\left(\left(\frac{-1-y}{\varepsilon}, \frac{1-y}{\varepsilon}\right)^c\right). \end{aligned}$$

Symmetry of the jump measure implies $\inf_{y \in \mathbb{R}} \nu_\eta^\varepsilon((\frac{-1-y}{\varepsilon}, \frac{1-y}{\varepsilon})^c) \geq \nu_\eta^\varepsilon([\frac{1}{\varepsilon}, \infty)) = \frac{1}{2}C_\varepsilon$. Using the elementary inequality $1 - e^{-x} \geq x - x^2/2$, $x \geq 0$, yields

$$(4.4) \qquad \mathbf{P}(\sigma_x < T) \geq T\left(1 - \frac{\beta_\varepsilon T}{2}\right)\frac{C_\varepsilon}{2}.$$

Then we apply Lemma 3.1(1) to obtain for any $t \geq 0$, $T \in (0, 2/\beta_\varepsilon)$ that

$$(4.5) \qquad \sup_{x \in I} \mathbf{P}(\sigma_x > t) \leq \left(1 - T\left(1 - \frac{\beta_\varepsilon T}{2}\right)\frac{C_\varepsilon}{2}\right)^{-1} \exp\left[-t\left(1 - \frac{\beta_\varepsilon T}{2}\right)\frac{C_\varepsilon}{2}\right].$$

The upper bound in Theorem 2.1 follows immediately by taking the infimum of the right-hand side of the latter estimate over $T \in (0, 2/\beta_\varepsilon)$.

4.2. *Super-exponential tails. Proof of Theorem 2.2, upper bound.* We need to show that for any $\delta \in (0, 1]$ there exists $T > 0$ and $\varepsilon_0 > 0$, such that for any $0 < \varepsilon < \varepsilon_0$ the following estimate holds:

$$(4.6) \qquad \inf_{x \in I} \mathbf{P}(\sigma_x \leq T) \geq TD_\varepsilon^{1+\delta}$$

with $D_\varepsilon$ defined in Theorem 2.2. Indeed, since $(1 - TD_\varepsilon^{1+\delta})^{-1} < 1 + \delta$ for $\varepsilon$ small enough, Lemma 3.1 yields the assertion.

Let $\delta' := \delta/7$, $M := \sup_{y \in I} |U'(y)|$, and let us choose $T \in (0, \frac{\delta'}{M} \wedge 1)$. Due to symmetry, it is sufficient to consider $x \in [0, 1)$. We start by remarking

that $\sigma_x > T$ implies the inequality

$$(4.7) \quad 1 > X_T^\varepsilon(x) = x - \int_0^T U'(X_{s-}^\varepsilon(x))\,ds + \varepsilon L_T \geq -MT + \varepsilon L_T.$$

Since $MT < \delta'$, we conclude $\{\sigma_x > T\} \subseteq \{\varepsilon L_T < 1 + \delta'\}$ and, thus,

$$(4.8) \quad \{\sigma_x \leq T\} \supseteq \{\varepsilon L_T \geq 1 + \delta'\}.$$

Recall that $q_\varepsilon$ denotes the $\varepsilon$-quantile of $\nu([\cdot, \infty))$, the positive tail of the Lévy measure, and set

$$(4.9) \quad g_\varepsilon := (\alpha - 1)^{-1/\alpha} q_\varepsilon.$$

Since the exponent $f(u) := -\ln \nu([u, \infty))$ is regularly varying at $+\infty$ with index $\alpha > 1$, we have $g_\varepsilon \to \infty$ and $\varepsilon g_\varepsilon \to 0$ as $\varepsilon \to 0$.

Consider a decomposition $L = \xi + \varphi^\varepsilon + \eta^\varepsilon$, where $\xi$, $\varphi^\varepsilon$ and $\eta^\varepsilon$ are Lévy processes having generating triplets $(d, \nu|_{(-\infty,1]}, \mu)$, $(0, \nu|_{(1,g_\varepsilon)}, 0)$ and $(0, \nu|_{[g_\varepsilon,\infty)}, 0)$ accordingly. If $N^{(\eta)}$ denotes the counting process of $\eta^\varepsilon$ and $\beta_\varepsilon^{(\eta)} = \nu([g_\varepsilon, \infty))$, we have $\varphi_T^\varepsilon > 0$ and $\eta_T^\varepsilon \geq g_\varepsilon N_T^{(\eta)}$. Inequality (4.8) yields the inclusion

$$(4.10) \quad \begin{aligned} \{\sigma_x \leq T\} &\supseteq \{\varepsilon \xi_T \geq -\delta'\} \cap \{\varepsilon \eta_T^\varepsilon \geq 1 + 2\delta'\} \\ &\supseteq \{\varepsilon \xi_T \geq -\delta'\} \cap \left\{N_T^{(\eta)} \geq \frac{1 + 2\delta'}{\varepsilon g_\varepsilon}\right\}. \end{aligned}$$

The random variables $\xi_T$ and $N_T^{(\eta)}$ are independent and $\mathbf{P}(\varepsilon \xi_T \geq -\delta') \geq 1 - \delta'$ for $\varepsilon$ small enough. Let $k_\varepsilon := [\frac{1+2\delta'}{\varepsilon g_\varepsilon}] + 1$. In particular, this means that $k_\varepsilon/\beta_\varepsilon^{(\eta)} \to \infty$ as $\varepsilon \to 0$. Thus, with the help of the inequality $k! \leq \frac{1}{2} k^k$, $k \geq 2$, and the previous estimate, we get for $\varepsilon$ sufficiently small that

$$(4.11) \quad \begin{aligned} \mathbf{P}(\sigma_x \leq T) &\geq (1 - \delta') \mathbf{P}(N_T^{(\eta)} = k_\varepsilon) \\ &= (1 - \delta') \exp(-\beta_\varepsilon^{(\eta)} T) \frac{(\beta_\varepsilon^{(\eta)} T)^{k_\varepsilon}}{k_\varepsilon!} \\ &\geq \frac{1}{2} \frac{(\beta_\varepsilon^{(\eta)})^{(1+\delta')k_\varepsilon}}{k_\varepsilon!} \\ &\geq \exp(-(1+\delta')k_\varepsilon(\ln k_\varepsilon + |\ln \beta_\varepsilon^{(\eta)}|)). \end{aligned}$$

Moreover, using the definition of regularly varying functions, and the fact that $q_\varepsilon \to \infty$ as $\varepsilon \to 0$, we estimate for small $\varepsilon$

$$(4.12) \quad |\ln \beta_\varepsilon^{(\eta)}| = f(g_\varepsilon) \leq (1 + \delta') \left(\frac{g_\varepsilon}{q_\varepsilon}\right)^\alpha f(q_\varepsilon) \leq \frac{1 + \delta'}{\alpha - 1} |\ln \varepsilon|.$$

By definition of $k_\varepsilon$, we have $\ln k_\varepsilon \leq |\ln \varepsilon|$ and $k_\varepsilon \leq \frac{1+3\delta'}{\varepsilon g_\varepsilon}$ for $\varepsilon$ small enough. This leads to the final estimate

$$\mathbf{P}(\sigma_x \leq T) \geq \exp\left(-(1+\delta')^2 \frac{\alpha}{\alpha-1} k_\varepsilon |\ln \varepsilon|\right)$$
(4.13)
$$\geq \exp\left(-(1+\delta')^2(1+3\delta')d_\alpha \frac{|\ln \varepsilon|}{\varepsilon q_\varepsilon}\right) \geq TD_\varepsilon^{1+\delta}.$$

4.3. *Bounded sub-exponential tails. Proof of Theorem 2.3, upper bound.*
Here we proceed as in the case of super-exponential tails. Let $\alpha \in (0,1)$ and $\theta \in (\frac{1}{n}, \frac{1}{n-1}]$ for some $n \geq 1$. Then $\phi(\theta) = (n-1)\theta^\alpha + \vartheta^\alpha$, $\vartheta = 1 - (n-1)\theta$.

Let $\delta > 0$ be fixed and let $\delta'$ be positive and specified later. Let $M = \sup_{x \in I} |U'(x)|$ and let $T \in (0, \frac{\delta'}{M} \wedge 1)$.

Consider a decomposition $L^{\varepsilon,\theta} = \xi^{\varepsilon,\theta} + \eta^{\varepsilon,\theta}$ as in Section 3.3 with big jumps $W_i$ being distributed with the law $\beta_{\varepsilon,\theta}^{-1} \nu^{\varepsilon,\theta}|_{[-g_\varepsilon, g_\varepsilon]^c}$, $\beta_{\varepsilon,\theta} = 2\nu((g_\varepsilon, \frac{\theta}{\varepsilon}])$, and $g_\varepsilon = 1/\sqrt{\varepsilon}$.

For $x \in [0, 1)$, we obtain similarly to (4.8) that

$$\{\sigma_x < T\} \supseteq \{\varepsilon L_T^{\varepsilon,\theta} > 1 + \delta'\}$$
$$\supseteq \{\varepsilon \xi_T^{\varepsilon,\theta} \geq -\delta'\} \cap \left\{\sum_{i=1}^n \varepsilon W_i \geq 1 + 2\delta'\right\} \cap \{N_T = n\}$$
(4.14)
$$\supseteq \{\varepsilon \xi_T^{\varepsilon,\theta} \geq -\delta'\} \cap \bigcap_{i=1}^{n-1} \{\varepsilon W_i \geq \theta - \delta'\}$$
$$\cap \{\varepsilon W_n \geq \vartheta + (n+1)\delta'\} \cap \{N_T = n\}.$$

We will take into account that $\mathbf{P}(\varepsilon \xi_T^{\varepsilon,\theta} \geq -\delta') \geq 1/2$ for $\varepsilon$ small enough. To estimate $\mathbf{P}(\varepsilon W_i > \theta - \delta')$ and $\mathbf{P}(\varepsilon W_n > \vartheta + (n+1)\delta')$, the following inequalities will be used. For any $0 < c < 1$, we have from the definition of regularly varying functions that for $c_1 < (1-c)^{-\alpha} - 1$ and $u$ big enough

(4.15) $\quad f(u) - f((1-c)u) \geq (1 - (1+c_1)(1-c)^\alpha)f(u) \geq u^{\alpha/2}$

and, thus,

(4.16) $\quad \frac{\nu([(1-c)u, u])}{\nu([u, \infty))} \geq \frac{\nu([(1-c)u, \infty))}{\nu([u, \infty))} - 1 = e^{f(u)-f((1-c)u)} - 1 \geq 1.$

Using independence of $\xi_T$, $W_i$ and $N_T$, the estimates $e^{-\beta_{\varepsilon,\theta} T} \geq 1/2$ and (4.16), and choosing $\delta' < \frac{\theta - \vartheta}{n+2}$ and small enough such that

(4.17) $\quad \frac{f((\vartheta + (n+2)\delta')/\varepsilon)}{f(\theta/\varepsilon)} \leq (1+\delta')\left(\frac{\vartheta + (n+2)\delta'}{\theta}\right)^\alpha \leq \frac{\delta}{3} + \left(\frac{\vartheta}{\theta}\right)^\alpha,$

we get for small $\varepsilon$

$$
\begin{aligned}
\mathbf{P}(\sigma_x &< T) \\
&\geq \frac{1}{2}\mathbf{P}(N_T = n)(\mathbf{P}(\varepsilon W_1 \geq \theta - \delta'))^{n-1}\mathbf{P}(\varepsilon W_1 \geq \vartheta + (n+1)\delta') \\
&\geq \frac{1}{2}e^{-\beta_{\varepsilon,\theta}T}\frac{T^n}{n!} \cdot \nu\left(\left[\frac{\theta - \delta'}{\varepsilon}, \frac{\theta}{\varepsilon}\right]\right)^{n-1} \cdot \nu\left(\left[\frac{\vartheta + (n+1)\delta'}{\varepsilon}, \frac{\theta}{\varepsilon}\right]\right) \\
(4.18) \quad &\geq \frac{1}{2}e^{-\beta_{\varepsilon,\theta}T}\frac{T^n}{n!} \cdot \nu\left(\left[\frac{\theta}{\varepsilon}, \infty\right)\right)^{n-1} \cdot \nu\left(\left[\frac{\vartheta + (n+2)\delta'}{\varepsilon}, \infty\right)\right) \\
&\geq \frac{T^n}{4n!}\exp\left(-(n-1)f\left(\frac{\theta}{\varepsilon}\right) - f\left(\frac{\vartheta + (n+2)\delta'}{\varepsilon}\right)\right) \\
&\geq T\exp\left(-\left(1 + \frac{\delta}{2}\right)\left(n - 1 + \left(\frac{\vartheta}{\theta}\right)^\alpha\right)f\left(\frac{\theta}{\varepsilon}\right)\right) \\
&\geq T\exp\left(-(1+\delta)((n-1)\theta^\alpha + \vartheta^\alpha)f\left(\frac{1}{\varepsilon}\right)\right) = TC_{\varepsilon,\theta}^{1+\delta}.
\end{aligned}
$$

## 5. The lower bounds.

5.1. *General remarks and reduction of starting values.* We will use the second part of Lemma 3.1 to establish the lower bound estimates for Theorems 2.1, 2.2 and 2.3. Consider $\sigma_x^-$ and $\sigma_x^*$ as in Lemma 3.2. In the sub-exponential case in Section 5.2 we will show that, for some appropriately chosen $\delta_0 > 0$, $p \geq 1$, $T > 0$ and any $\delta \in (0, \delta_0)$, the following estimate holds for small $\varepsilon$:

$$(5.1) \qquad \sup_{|x| \leq \delta} \mathbf{P}(\sigma_x^* < T) \leq TC_\varepsilon^{1-p\delta}.$$

In Sections 5.4 and 5.3 we establish the analogous inequalities with $C_{\varepsilon,\theta}$ and $D_\varepsilon$ replacing $C_\varepsilon$ for the case of bounded and super-exponential jumps, respectively. The estimate (5.1) established, thus Lemma 3.1 yields

$$
\begin{aligned}
(5.2) \quad \inf_{|x| \leq \delta} \mathbf{P}(\sigma_x > t) &\geq \inf_{|x| \leq \delta} \mathbf{P}(\sigma_x^- > t) \\
&\geq (1 - TC_\varepsilon^{1-p\delta})\exp\left(\frac{\ln(1 - TC_\varepsilon^{1-p\delta})}{T}t\right) \\
&\geq (1 - \delta)\exp(-C_\varepsilon^{1-(p+1)\delta}t)
\end{aligned}
$$

for $\varepsilon$ sufficiently small and uniformly over $t \geq 0$.

It is left to get rid of the constraint $|x| \leq \delta$ on the initial value. This can be done with help of the inequality

$$(5.3) \qquad \inf_{x \in I_\delta^-} \mathbf{P}(\sigma_x > t) \geq (1 - \delta) \inf_{|x| \leq \delta} \mathbf{P}(\sigma_x > t),$$

which will be proven to hold for all $\delta \in (0, 1/2]$ and $\varepsilon$ small enough. Indeed, let $\delta' = \frac{\delta}{p+1} \leq \frac{\delta}{2}$. Thus, (5.2) and (5.3) yield

$$\inf_{x \in I_\delta^-} \mathbf{P}(\sigma_x > t) \geq \inf_{x \in I_{\delta'}^-} \mathbf{P}(\sigma_x > t) \geq (1 - \delta') \inf_{|x| \leq \delta'} \mathbf{P}(\sigma_x > t)$$
(5.4)
$$\geq (1 - \delta')^2 \exp(-C_\varepsilon^{1-(p+1)\delta'} t) \geq (1 - \delta) \exp(-C_\varepsilon^{1-\delta} t)$$

for any $t \geq 0$ and $\varepsilon$ sufficiently small. This entails the lower bounds of the estimate (2.7). The estimates leading to (2.9) and (2.16) are obtained analogously.

The structures of the proofs providing the lower bounds in Theorems 2.1, 2.2 and 2.3 are similar. As is shown above, it is sufficient to prove that inequalities (5.1) and (5.3) hold. To do this, we consider a decomposition $L = \xi^\varepsilon + \eta^\varepsilon$ as in Section 3.3 with some appropriately chosen $g_\varepsilon$, such that $g_\varepsilon \to \infty$ and $\varepsilon g_\varepsilon \to 0$ as $\varepsilon \to 0$. In the Sections 5.2.1, respectively, 5.3.1 we will use Lemma 3.2 and Corollary 3.1 to obtain an embedding of $\{\sigma_x^* < T\}$ in terms of sets described by the large and small jump parts $\xi^\varepsilon$ and $\eta^\varepsilon$. In Sections 5.2.2, respectively, 5.3.2 we will use Lemmas 3.3 and 3.4 to estimate the probabilities of the covering sets.

PROOF OF INEQUALITY (5.3). Let $\sigma_x^1 := \inf\{t \geq 0 : |X_t^\varepsilon(x)| \leq \delta\}$ and $\sigma_x^2 := \inf\{t \geq 0 : |X_t^\varepsilon(x)| \geq 1 - \frac{\delta}{2}\}$. The strong Markov property and time homogeneity of $X^\varepsilon$ yield for any $x \in I_\delta^-$ that

$$\mathbf{P}(\sigma_x > t) \geq \mathbf{P}(\sigma_x > \sigma_x^1 + t) \geq \mathbf{P}(\sigma_x > \sigma_x^1) \inf_{|x| \leq \delta} \mathbf{P}(\sigma_x > t)$$
(5.5)
$$\geq \mathbf{P}(\sigma_x^2 > \sigma_x^1) \inf_{|x| \leq \delta} \mathbf{P}(\sigma_x > t).$$

Let $\tilde{m} := \min_{y \in I_{\delta/2}^- \setminus [-\delta, \delta]} |U'(y)|$ and $T := 2/\tilde{m}$. We have

(5.6) $\{\sigma_x^2 \leq \sigma_x^1\} \subseteq \{\sigma_x^1 \geq T, \sigma_x^2 \geq T, \tau_1 > T\} \cup \{\sigma_x^2 < T \wedge \tau_1\} \cup \{\tau_1 \leq T\}.$

We choose $\varepsilon$ sufficiently small such that $\varepsilon g_\varepsilon < \delta$. For such $\varepsilon$ in analogy with (3.9), we get $\{\sigma_x^1 \geq T, \sigma_x^2 \geq T, \tau_1 > T\} \subseteq \{\sup_{t \leq T} |\varepsilon \xi_t^\varepsilon| > \tilde{m}T - 1 = 1\}$, and for any $x \in I_\delta^-$, the second part of Lemma 3.2 yields $\{\sigma_x^2 < T \wedge \tau_1\} \subseteq \{(\sup - \inf_{t < T \wedge \tau_1}) \varepsilon L_t \geq \frac{\delta}{2}\} \subseteq \{\sup_{t \leq T} |\varepsilon \xi_t^\varepsilon| \geq \frac{\delta}{4}\}$. Thus, since $g_\varepsilon \to \infty$ and $\varepsilon g_\varepsilon \to 0$ as $\varepsilon \to 0$, we have

(5.7) $$\sup_{x \in I_\delta^-} \mathbf{P}(\sigma_x^2 \leq \sigma_x^1) \leq \mathbf{P}\left(\sup_{t \leq T} |\varepsilon \xi_t^\varepsilon| \geq \frac{\delta}{4}\right) + \mathbf{P}(\tau_1 \leq T) \leq \delta$$

for $\varepsilon$ small enough, and (5.5) yields the assertion. $\square$

5.2. *Sub-exponential tails. Proof of Theorem 2.2, lower bound.*

5.2.1. *Estimate of beginning of exit.* Let $J_\delta = [-\delta, \delta]$, and $m = \inf\{|U'(y)| : y \in I_\delta^- \setminus J_\delta\}$ as in Lemma 3.2. The following estimates would hold analogously for any choice $T > 0$ and $\hat{T} > \frac{1}{m}$. For simplicity, we choose $\hat{T} = \frac{2}{m}$ and $T = \hat{T}$, such that $m\hat{T} - 1 = 1$ and $T + \hat{T} = 2T$. Let $k \geq 1$ and abbreviate

$$\chi_k := \bigcup_{i=0}^{k-1} \left\{ \sup_{t \leq T} \varepsilon |\xi^\varepsilon_{S_i + T + t} - \xi^\varepsilon_{S_i + T}| \geq 1 \right\}. \tag{5.8}$$

For $x \in J_\delta$, we have

$$\{\sigma^*_x < T\} \subseteq \chi_k \cup (\{\sigma^*_x < T, \sigma^-_x \geq S_k \wedge 2Tk\} \setminus \chi_k) \tag{5.9}$$
$$\cup \bigcup_{j=0}^{k-1} (\{\sigma^*_x < T, S_j \leq \sigma^-_x < S_{j+1} \wedge 2Tk\} \setminus \chi_k).$$

For $1 \leq j \leq k$, the first part of Corollary 3.1 yields that

$$\{\sigma^* < T, \sigma^- \geq S_j\} \subseteq \chi_j \cup \bigcap_{i=1}^{j} \{\tau_i \leq 2T\} \tag{5.10}$$

and for all $1 \leq j \leq k-1$ and $x \in J_\delta$, we obtain with the help of Corollary 3.1(ii) and the previous inclusion that

$$\{\sigma^*_x < T, S_j \leq \sigma^-_x < S_{j+1} \wedge 2Tk\} \setminus \chi_k$$
$$\subseteq \{\sigma^-_x < S_{j+1} \wedge 2Tk\} \cap (\{\sigma^*_x < T, \sigma^-_x \geq S_j\} \setminus \chi_j) \tag{5.11}$$
$$\subseteq \left\{ \sup_{t < S_{j+1} \wedge 2Tk} |X^\varepsilon_t(x)| \geq 1 - \delta \right\} \cap \bigcap_{i=1}^{j} \{\tau_i \leq 2T\}$$
$$\subseteq \left( \left\{ \sum_{i=1}^{j} |\varepsilon W_i| \geq 1 - 3\delta \right\} \cup \left\{ \sup_{t \leq 2Tk} |\varepsilon \xi^\varepsilon_t| \geq \frac{\delta}{2} \right\} \right) \cap \bigcap_{i=1}^{j} \{\tau_i \leq 2T\}.$$

In the particular case $j = 0$ we get directly with the help of Corollary 3.1(ii) that

$$\{\sigma^*_x < T, \sigma^-_x \in [0, S_1 \wedge 2Tk)\} \setminus \chi_k \subseteq \{\sigma^*_x < T, \sigma^-_x < S_1 \wedge 2Tk\} \tag{5.12}$$
$$\subseteq \left\{ \sup_{t \leq 2Tk} |\varepsilon \xi^\varepsilon_t| \geq \frac{\delta}{2} \right\}.$$

Putting all sets in (5.9) together, we obtain

$$\{\sigma^*_x < T\} \subseteq \chi_k \cup \bigcap_{i=1}^{k} \{\tau_i \leq 2T\} \cup \left\{ \sup_{t \leq 2kT} |\varepsilon \xi^\varepsilon_t| \geq \frac{\delta}{2} \right\} \tag{5.13}$$
$$\cup \bigcup_{j=1}^{k-1} \left( \bigcap_{i=1}^{j} \{\tau_i \leq 2T\} \cap \left\{ \sum_{i=1}^{j} |\varepsilon W_i| \geq 1 - 3\delta \right\} \right).$$

5.2.2. *The lower bound via beginning of exit.* Let $\delta_0 := \frac{1}{3}\alpha(1-\alpha)$. Pick $\delta \in (0, \delta_0)$ and let

(5.14) $$g_\varepsilon := \varepsilon^{-(1-\alpha-\delta)} \quad \text{and} \quad k = k_\varepsilon := [\varepsilon^{-(\alpha-\delta)}].$$

Then (5.13) yields

(5.15) $$\mathbf{P}(\sigma_x^* < T) \leq \mathbf{P}(\chi_{k_\varepsilon}) + \mathbf{P}(\tau_1 \leq 2T)^{k_\varepsilon} + \mathbf{P}\left(\sup_{t \leq 2k_\varepsilon T} |\varepsilon \xi_t^\varepsilon| \geq \frac{\delta}{2}\right)$$
$$+ \sum_{k=1}^{k_\varepsilon - 1} \mathbf{P}(\tau_1 \leq 2T)^k \cdot \mathbf{P}\left(\sum_{i=1}^k |\varepsilon W_i| \geq 1 - 3\delta\right).$$

In the next steps we estimate the summands of the previous formula. Recall that by definition we have $\varepsilon k_\varepsilon g_\varepsilon \to 0$ as $\varepsilon \to 0$.

1. We first apply the strong Markov property of $\xi^\varepsilon$ and Lemma 3.3 with $2k_\varepsilon T$ instead of $T$ and $\frac{\delta}{2\varepsilon}$ instead of $f$ to the first and third term on the right-hand side of (5.15) to obtain for $\varepsilon$ sufficiently small

(5.16)
$$\mathbf{P}(\chi_{k_\varepsilon}) + \mathbf{P}\left(\sup_{t \leq 2k_\varepsilon T} |\varepsilon \xi_t^\varepsilon| \geq \frac{\delta}{2}\right)$$
$$\leq k_\varepsilon \mathbf{P}\left(\sup_{t \leq T} |\varepsilon \xi_t^\varepsilon| \geq 1\right) + \mathbf{P}\left(\sup_{t \leq 2k_\varepsilon T} |\varepsilon \xi_t^\varepsilon| \geq \frac{\delta}{2}\right)$$
$$\leq (k_\varepsilon + 1) \mathbf{P}\left(\sup_{t \leq 2k_\varepsilon T} |\varepsilon \xi_t^\varepsilon| \geq \frac{\delta}{2}\right)$$
$$\leq (k_\varepsilon + 1) \exp\left(-(1-\delta)\frac{\delta}{2\varepsilon g_\varepsilon} \ln \frac{\delta}{4T\varepsilon k_\varepsilon g_\varepsilon}\right)$$
$$\leq \exp\left(-\frac{1}{\varepsilon g_\varepsilon} |\ln \varepsilon|\right)$$
$$\leq \exp(-\varepsilon^{-(\alpha+\delta/2)}) \leq TC_\varepsilon.$$

2. We next deal with the second term on the right-hand side of (5.15). Recall that the negative logarithm of the Lévy measure's tail $f(u)$, $u > 1$, is a regularly varying function with index $\alpha \in (0, 1)$. Since $g_\varepsilon \to \infty$ as $\varepsilon \to 0$, we can choose $\varepsilon$ sufficiently small, such that $f(g_\varepsilon) \geq g_\varepsilon^{\alpha-\delta}$ and the following estimate holds:

(5.17)
$$\mathbf{P}(\tau_1 \leq 2T)^{k_\varepsilon} \leq (2T\beta_\varepsilon)^{k_\varepsilon} = (4Te^{-f(g_\varepsilon)})^{k_\varepsilon} \leq (4Te^{-g_\varepsilon^{\alpha-\delta}})^{k_\varepsilon}$$
$$\leq \exp(-(1-\delta)\varepsilon^{-(1-\alpha-\delta)(\alpha-\delta)-(\alpha-\delta)})$$
$$\leq \exp(-(1-\delta)\varepsilon^{-(\alpha+\alpha(1-\alpha)-3\delta)}).$$

For the first inequality in the chain (5.17), we hereby use that the law of $\tau_1$ is exponential with mean $\beta_\varepsilon^{-1}$, while $\alpha \in (0,1)$ is needed for the last. Thus, the hypothesis $\delta < \frac{1}{3}\alpha(1-\alpha)$ yields $\mathbf{P}(\tau_1 \leq 2T)^{k_\varepsilon} < TC_\varepsilon$ for $\varepsilon$ sufficiently small.

3. We finally treat a summand of the last term on the right-hand side of (5.15). By definition, we have $\varepsilon g_\varepsilon k_\varepsilon \to 0$ as $\varepsilon \to 0$. Thus, $\varepsilon$ can be chosen sufficiently small, such that by Lemma 3.4 applied with $r = \frac{1-3\delta}{\varepsilon}$, to estimate $\mathbf{P}(\sum_{i=1}^k |\varepsilon W_i| \geq 1 - 3\delta)$, we get the following inequalities:

$$
\begin{aligned}
&\max_{1 \leq k \leq k_\varepsilon - 1} \mathbf{P}(\tau_1 \leq 2T)^k \cdot \mathbf{P}\bigg(\sum_{i=1}^k |\varepsilon W_i| \geq 1 - 3\delta\bigg) \\
&\leq \max_{1 \leq k \leq k_\varepsilon - 1} (2T\beta_\varepsilon)^k \cdot \mathbf{P}\bigg(\sum_{i=1}^k |\varepsilon W_i| \geq 1 - 3\delta\bigg) \\
&\leq \max_{1 \leq k \leq k_\varepsilon - 1} \bigg(8T + 4T\frac{|\ln \varepsilon| - \ln g_\varepsilon}{\ln(1+\delta)}\bigg)^k \\
&\quad \times \exp\bigg(-\inf\bigg\{\sum_{i=1}^k f(x_i) : \sum_{i=1}^k x_i = \frac{(1-3\delta)(1-\delta)}{\varepsilon}, \\
&\hspace{6cm} x_i \in \bigg[g_\varepsilon, \frac{1-3\delta}{\varepsilon}\bigg]\bigg\}\bigg) \\
&\leq |\ln \varepsilon|^{2k_\varepsilon} \\
&\quad \times \max_{1 \leq k \leq k_\varepsilon - 1} \exp\bigg(-\inf\bigg\{\sum_{i=1}^k f(x_i) : \sum_{i=1}^k x_i = \frac{(1-3\delta)(1-\delta)}{\varepsilon}, \\
&\hspace{6cm} x_i \in \bigg[g_\varepsilon, \frac{1}{\varepsilon}\bigg]\bigg\}\bigg).
\end{aligned}
$$
(5.18)

Note that for the last inequality, the minimizer on the intervals $[g_\varepsilon, \frac{1}{\varepsilon}]$ is smaller, so it gives an upper estimate. For the crucial estimate of the exponential rate, we call upon Potter's bound for the regularly varying function $f$ (see [2], Theorem 1.5.6). Since $g_\varepsilon \to \infty$ as $\varepsilon \to 0$, it provides the following estimate for all $x \in [g_\varepsilon, \frac{1}{\varepsilon}]$ and $\varepsilon$ sufficiently small:

(5.19) $$f(x) \geq (1-\delta)f(\varepsilon^{-1})(\varepsilon x)^{\alpha+\delta}.$$

Therefore, we get

$$
\begin{aligned}
&\inf\bigg\{\sum_{i=1}^k f(x_i) : \sum_{i=1}^k x_i = (1-3\delta)(1-\delta)\varepsilon^{-1}, x_i \in [g_\varepsilon, \varepsilon^{-1}]\bigg\} \\
&\geq (1-\delta)f(\varepsilon^{-1})\varepsilon^{\alpha+\delta} \\
&\quad \times \inf\bigg\{\sum_{i=1}^k x_i^{\alpha+\delta} : \sum_{i=1}^k x_i = (1-3\delta)(1-\delta)\varepsilon^{-1}, x_i \in [g_\varepsilon, \varepsilon^{-1}]\bigg\} \\
&\geq (1-\delta)f(\varepsilon^{-1})\varepsilon^{\alpha+\delta}((1-3\delta)(1-\delta)\varepsilon^{-1})^{\alpha+\delta} \geq (1-5\delta)f(\varepsilon^{-1}).
\end{aligned}
$$
(5.20)

For obtaining the second inequality in the preceding chain, we have to recall that $\alpha + \delta < 1$, and use the inequality (3.34). So the estimation in (5.18) may be completed by

$$
\begin{aligned}
\max_{1 \leq k \leq k_\varepsilon - 1} &\mathbf{P}(\tau_1 \leq 2T)^k \cdot \mathbf{P}\left(\sum_{i=1}^k |\varepsilon W_i| \geq 1 - 3\delta\right) \\
&\leq \exp(-(1-5\delta)f(\varepsilon^{-1}) + 2k_\varepsilon \ln|\ln \varepsilon|) \\
&\leq \exp(-(1-6\delta)f(\varepsilon^{-1})) \leq TC_\varepsilon^{1-7\delta},
\end{aligned}
\tag{5.21}
$$

which again holds for $\varepsilon$ small enough. Collecting the bounds we obtained for the terms in (5.15), we finally get for $\varepsilon$ small enough

$$
\mathbf{P}(\sigma_x^* < T) \leq 2TC_\varepsilon + (k_\varepsilon - 1)TC_\varepsilon^{1-7\delta} \leq TC_\varepsilon^{1-8\delta}.
\tag{5.22}
$$

Thus, we obtain the desired upper bound, and complete the proof of Theorem 2.1.

### 5.3. Super-exponential tails. Proof of Theorem 2.2, lower bound.

#### 5.3.1. Estimate of beginning of exit.

Again, we start by covering the crucial set $\{\sigma_x^* < T\}$, $x \in J_\delta$, by sets described in terms of the small and large jump parts. Let $k \geq 1$. Let $m$, $\hat{T}$, $T$ and $\chi_k$ be defined as in Section 5.2.1. Then with help of Corollary 3.1(i), we obtain

$$
\begin{aligned}
\{\sigma_x^* < T\} &= \{\sigma_x^* < T, \sigma_x^- < 2Tk \wedge S_k\} \cup \{\sigma_x^* < T, \sigma_x^- \geq 2Tk \wedge S_k\} \\
&\subseteq \{\sigma_x^* < T, \sigma_x^- < 2Tk \wedge S_k\} \cup \chi_k \cup \bigcap_{i=1}^k \{\tau_i \leq 2T\}.
\end{aligned}
\tag{5.23}
$$

Define the set

$$
\bar{\chi}_k := \bigcup_{i=1}^{k-1} \{\varepsilon W_i \geq \delta\}.
\tag{5.24}
$$

We notice that on the event $\{\sigma_x^* < T, \sigma_x^- < 2Tk \wedge S_k\} \cap \bar{\chi}_k^c$ the estimate

$$
\left(\sup_{t \leq \sigma_x^-} - \inf\right) \varepsilon L_t \geq m(\sigma_x^- - \sigma_x^*)
\tag{5.25}
$$

holds. To see this, we use similar arguments as in the proof of Lemma 3.2(i). Indeed, on the event $\{\sigma_x^* < T, \sigma_x^- < 2Tk \wedge S_k\} \cap \bar{\chi}_k^c$ the process $X^\varepsilon(x)$ does not change its sign during the time interval $[\sigma_x^*, \sigma_x^-]$, and for any $t \in (\sigma_x^*, \sigma_x^-)$ we have $|U'(X_{t-}(x))| \geq m$. Further, we have that if $0 < X_{\sigma_x^*-}(x) \leq \delta$, then $X_{\sigma_x^-}(x) \geq 1-\delta$, and $X_{\sigma_x^-}(x) = X_{\sigma_x^*-}(x) + \int_{\sigma_x^*}^{\sigma_x^-} U'(X_{t-}(x))\,dt + \varepsilon(L_{\sigma_x^-} - L_{\sigma_x^*-})$, and, thus,

$$(5.26) \quad \left(\sup_{t \leq \sigma_x^-} - \inf\right)\varepsilon L_t \geq (1-2\delta) + m(\sigma_x^- - \sigma_x^*) \geq m(\sigma_x^- - \sigma_x^*).$$

The case of negative values $X_{\sigma_x^*-}(x)$ is considered analogously.

This and Lemma 3.2(ii) lead to the following estimate for $k \geq 1$ and $x \in J_\delta$ (recall $mT = 2$):

$$(5.27) \quad \begin{aligned}
&\{\sigma_x^* < T, \sigma_x^- < 2Tk \wedge S_k\} \cap \bar\chi_k^c \\
&\subseteq \{\sigma_x^- < 4T\} \\
&\quad \cup \bigcup_{i=2}^{k-1}(\{\sigma_x^* < T, \sigma_x^- \in [2iT, 2(i+1)T], \sigma_x^- < S_k\} \cap \bar\chi_k^c) \\
&\subseteq \left\{\left(\sup_{t<4T} - \inf\right)\varepsilon L_t \geq 1-2\delta\right\} \\
&\quad \cup \bigcup_{i=2}^{k-1}\left\{\left(\sup_{t<2(i+1)T} - \inf\right)\varepsilon L_t \geq (2i-1)Tm\right\} \\
&\subseteq \left\{\left(\sup_{t<4T} - \inf\right)\varepsilon L_t \geq 1-2\delta\right\} \\
&\quad \cup \bigcup_{i=2}^{k-1}\left\{\sum_{j=0}^{i}\left(\sup_{t<2T} - \inf\right)\varepsilon(L_{2jT+t} - L_{2jT}) \geq i+1\right\} \\
&\subseteq \bigcup_{i=0}^{k-1}\left\{\left(\sup_{t<4T} - \inf\right)\varepsilon(L_{2iT+t} - L_{2iT}) \geq 1-2\delta\right\}.
\end{aligned}$$

In particular, this entails

$$(5.28) \quad \begin{aligned}
\{\sigma_x^* < T\} &\subseteq \chi_k \cup \bar\chi_k \cup \bigcap_{i=1}^{k}\{\tau_i \leq 2T\} \\
&\quad \cup \bigcup_{i=0}^{k-1}\left\{\left(\sup_{t \leq 4T} - \inf\right)\varepsilon(L_{2iT+t} - L_{2iT}) \geq 1-2\delta\right\}.
\end{aligned}$$

5.3.2. *The lower bound via beginning of exit.* Let $\delta \in (0, \frac{1}{6} \wedge (\alpha-1))$ be fixed, let $q_\varepsilon := \sup\{u > 0 : \nu([u, \infty)) \geq \varepsilon\}$ denote the $\varepsilon$-quantile of the Lévy measure $\nu$, and set

$$(5.29) \quad g_\varepsilon := \frac{q_\varepsilon}{3} \quad \text{and} \quad k_\varepsilon := \left[\frac{|\ln \varepsilon|}{\varepsilon}\right].$$

In particular, since the tails of $\nu$ are super-exponential, $\varepsilon g_\varepsilon \to 0$ as $\varepsilon \to 0$.

We shall also use a simple estimate of $d_\alpha := \inf_{y>0}(y^{-\alpha+1} + y) = \alpha(\alpha - 1)^{1/\alpha-1}$. It is easy to see that

$$(5.30) \quad 1 = \inf_{y>0}(y^{-\alpha+1} \vee y) < d_\alpha \leq (y^{-\alpha+1} + y)|_{y=1} = 2.$$

Next, we estimate the probabilities of the events in (5.28).

1. We use the inequality $d_\alpha \leq 2 \leq 3(1 - 2\delta)$ and Lemma 3.3 with $f = \frac{1}{\varepsilon}$ to obtain for $\varepsilon$ small that

$$(5.31) \quad \mathbf{P}(\chi_k) \leq k_\varepsilon \mathbf{P}\left(\sup_{t \leq 2T} |\varepsilon \xi_t^\varepsilon| \geq 1\right) \leq k_\varepsilon \exp\left(-\frac{1-\delta}{\varepsilon g_\varepsilon} \ln \frac{1}{2T\varepsilon g_\varepsilon}\right)$$

$$\leq k_\varepsilon \exp(-3(1-2\delta)\varepsilon^{-1} q_\varepsilon^{-1} |\ln \varepsilon|) < TD_\varepsilon.$$

2. To deal with the second term, recall that $f(u) = -\ln \nu([u, \infty))$ is regularly varying with $\alpha > 1$ at infinity, $g_\varepsilon < q_\varepsilon$ and, thus, $\beta_\varepsilon > \varepsilon$. Then we have

$$(5.32) \quad \mathbf{P}(\bar{\chi}_{k_\varepsilon}) \leq k_\varepsilon \mathbf{P}(|\varepsilon W_1| \geq \delta) \leq 2k_\varepsilon \beta_\varepsilon^{-1} e^{-f(\delta/\varepsilon)}$$

$$\leq \exp(-(\delta/\varepsilon)^{\alpha-\delta} + 2|\ln \varepsilon| + \ln k_\varepsilon) \leq TD_\varepsilon.$$

3. Since $\beta_\varepsilon \to 0$ as $\varepsilon \to 0$, we have for sufficiently small $\varepsilon > 0$

$$(5.33) \quad \mathbf{P}\left(\bigcap_{i=1}^{k_\varepsilon}\{\tau_i \leq 2T\}\right) = \mathbf{P}(\tau_1 \leq 2T)^{k_\varepsilon} \leq (2T\beta_\varepsilon)^{k_\varepsilon} \leq \exp(-k_\varepsilon) < TD_\varepsilon.$$

4. The estimate for the last union in (5.28) is the most important part of the proof. Since for any $0 \leq i \leq k-1$ the processes $(L_{2iT+t} - L_{2iT})_{t\geq 0}$ have the same law as $L = (L_t)_{t\geq 0}$, its enough to work with the original process $L$. We prove the following lemma.

LEMMA 5.1. *Let the jump measure $\nu$ of $L$ be symmetric and super-exponential with index $\alpha > 1$. Then for any $T > 0$, $a > 0$, and $\delta > 0$, there exists $\varepsilon_0 > 0$ such that for all $0 < \varepsilon \leq \varepsilon_0$ the following estimate holds:*

$$(5.34) \quad \mathbf{P}\left(\left(\sup_{t \leq T} - \inf_{t \leq T}\right)\varepsilon L_t > a\right) \leq D_\varepsilon^{(1-\delta)a}.$$

Applying Lemma 5.1 with $a = 1 - 2\delta$ to the last event in (5.28), we get

$$(5.35) \quad k_\varepsilon \mathbf{P}\left(\left(\sup_{t \leq 4T} - \inf\right)\varepsilon L_t \geq 1 - 2\delta\right) \leq k_\varepsilon D_\varepsilon^{1-3\delta} \leq TD_\varepsilon^{1-4\delta}$$

for sufficiently small $\varepsilon$. This completes the proof of Theorem 2.2.

PROOF OF LEMMA 5.1. Due to monotonicity, it is sufficient to consider $\delta \in (0, \frac{1}{21\alpha} \wedge a)$. Let $p := 9\alpha$. We shall prove that for any such $\delta$ there exists $\varepsilon_0 > 0$, such that for every $0 < \varepsilon < \varepsilon_0$ the estimate $\mathbf{P}((\sup - \inf_{t \leq T})\varepsilon L_t > a) \leq D_\varepsilon^{(1-p\delta)a}$ holds. This entails the asserted inequality.

Consider a decomposition $L = \xi^\varepsilon + \eta^\varepsilon$ as in Section 3.3 with the threshold

$$g_\varepsilon := \frac{\delta q_\varepsilon}{6}. \tag{5.36}$$

Note that this $g_\varepsilon$ is different from its counterpart defined in (5.29) at the beginning of this subsection, and will be only used in the proof of the lemma.

Since $(\sup - \inf_{t \leq T})\varepsilon \xi_t^\varepsilon \leq 2 \sup_{t \leq T} |\varepsilon \xi_t^\varepsilon|$, we have for any $n \geq 1$

$$\begin{aligned}
&\left\{\left(\sup_{t \leq T} - \inf\right) \varepsilon L_t \geq a\right\} \\
&\subseteq \left\{\left(\sup_{t \leq T} - \inf\right) \varepsilon \xi_t^\varepsilon \geq \delta a\right\} \cup \left\{\left(\sup_{t \leq T} - \inf\right) \varepsilon \eta_t^\varepsilon \geq (1-\delta)a\right\} \\
&\subseteq \left\{\sup_{t \leq T} |\varepsilon \xi_t^\varepsilon| \geq \frac{1}{2}\delta a\right\} \cup \{N_T > n\} \\
&\quad \cup \bigcup_{k=1}^{n}\left\{N_T = k, \sum_{i=1}^{k} |\varepsilon W_i| > (1-\delta)a\right\}.
\end{aligned} \tag{5.37}$$

The goal of the next steps consists in estimating the probabilities of the events figuring in the second line of (5.37) with an appropriately chosen $n$, namely, with

$$n := n_\varepsilon = \left[\frac{3a}{q_\varepsilon \varepsilon}\right]. \tag{5.38}$$

1. Recall that $d_\alpha \leq 2$ and apply Lemma 3.3 with $f = \frac{\delta a}{2\varepsilon}$ to get

$$\begin{aligned}
\mathbf{P}\left(\sup_{t \leq T} |\varepsilon \xi_t^\varepsilon| \geq \frac{\delta a}{2}\right) &\leq \exp\left(-(1-\delta)\frac{3a}{\varepsilon q_\varepsilon} \ln \frac{3a}{T \varepsilon q_\varepsilon}\right) \\
&\leq \exp\left(-(1-2\delta)\frac{3a}{\varepsilon q_\varepsilon} |\ln \varepsilon|\right) < D_\varepsilon^a
\end{aligned} \tag{5.39}$$

for $\varepsilon$ small enough.

2. To estimate $\mathbf{P}(N_T > n_\varepsilon)$, we will use Stirling's formula. By choice of $g_\varepsilon$ and $n_\varepsilon$ we have $\varepsilon g_\varepsilon \to 0$, $\frac{\ln n_\varepsilon}{|\ln \varepsilon|} \to 1$, and $\beta_\varepsilon T \leq 1$ as $\varepsilon \to 0$. Thus, for $\varepsilon$ sufficiently small, the estimate $n_\varepsilon! \geq \exp(n_\varepsilon(\ln n_\varepsilon - 1)) \geq \exp((1-\delta)n_\varepsilon|\ln \varepsilon|)$ holds, and we get

$$\begin{aligned}
\mathbf{P}(N_T > n_\varepsilon) &\leq \sum_{k=n_\varepsilon}^{\infty} \frac{(\beta_\varepsilon T)^k}{k!} \leq \frac{(\beta_\varepsilon T)^{n_\varepsilon}}{n_\varepsilon!} \sum_{k=0}^{\infty} \left(\frac{\beta_\varepsilon T}{n_\varepsilon}\right)^k \leq (1+\delta)\frac{(\beta_\varepsilon T)^{n_\varepsilon}}{n_\varepsilon!} \\
&\leq \exp(-(1-\delta)n_\varepsilon|\ln \varepsilon|)
\end{aligned} \tag{5.40}$$

$$\leq \exp\left(-(1-2\delta)\frac{3a}{\varepsilon q_\varepsilon}|\ln\varepsilon|\right) < D_\varepsilon^a.$$

3. As in Section 5.2, the crucial ingredient which produces the phase transition at $\alpha = 1$ comes from the following estimate for the exponential rate of sums of large jumps. Recall that $f(u) = -\ln\nu([u,+\infty))$, $u > 0$, is regularly varying with index $\alpha > 1$. By choice of parameters, we have $\varepsilon g_\varepsilon n_\varepsilon \to \frac{\delta}{2} < (1-\delta)^2 a$ as $\varepsilon \to 0$. Thus, Lemma 3.4 with $r = \frac{(1-\delta)a}{\varepsilon}$ can be used to estimate $\mathbf{P}(\sum_{i=1}^k |\varepsilon W_i| \geq (1-\delta)a)$ for $1 \leq k \leq n_\varepsilon - 1$. Hence, for $\varepsilon$ sufficiently small, the following estimate holds uniformly for all $1 \leq k < n_\varepsilon$:

$$\mathbf{P}\left(\sum_{i=1}^k |\varepsilon W_i| \geq (1-\delta)a\right)$$

(5.41)
$$\leq \beta_\varepsilon^{-k} |\ln\varepsilon|^{2k} \exp\left(-\inf\left\{\sum_{i=1}^k f(x_i) : \sum_{i=1}^k x_i = \frac{(1-\delta)^2 a}{\varepsilon}, \right.\right.$$
$$\left.\left. x_i \in \left[g_\varepsilon, \frac{(1-\delta)a}{\varepsilon}\right]\right\}\right).$$

Again we invoke Potter's bound to estimate the negative of the exponential rate. Choose $\tilde{\delta}$ sufficiently small, such that $\alpha - \tilde{\delta} > 1$ and $(\frac{\delta}{2})^{\tilde{\delta}} > 1 - \delta$. For sufficiently small $\varepsilon$, the estimate $f(x) \geq (1-\delta)f(g_\varepsilon)(\frac{x}{g_\varepsilon})^{\alpha-\tilde{\delta}}$ then holds for any $x \geq g_\varepsilon$. Thus, for any $1 \leq k < n_\varepsilon$, we get

$$\inf\left\{\sum_{i=1}^k f(x_i) : \sum_{i=1}^k x_i = \frac{(1-\delta)^2 a}{\varepsilon}, x_i \in \left[g_\varepsilon, \frac{(1-\delta)a}{\varepsilon}\right]\right\}$$

$$\geq (1-\delta)\frac{f(g_\varepsilon)}{g_\varepsilon^{\alpha-\tilde{\delta}}} \inf\left\{\sum_{i=1}^k x_i^{\alpha-\tilde{\delta}} : \sum_{i=1}^k x_i = \frac{(1-\delta)^2 a}{\varepsilon}, x_i > 0\right\}$$

$$= (1-\delta)\frac{f(g_\varepsilon)}{g_\varepsilon^{\alpha-\tilde{\delta}}} k \left(\frac{(1-\delta)^2 a}{\varepsilon k}\right)^{\alpha-\tilde{\delta}}$$

(5.42)
$$\geq (1-3\alpha\delta)f(g_\varepsilon) k \left(\frac{a}{\varepsilon g_\varepsilon k}\right)^{\alpha-\tilde{\delta}}$$

$$\geq (1-3\alpha\delta)f(g_\varepsilon) n_\varepsilon \left(\frac{a}{\varepsilon g_\varepsilon n_\varepsilon}\right)^{\alpha-\tilde{\delta}} \left(\frac{n_\varepsilon}{k}\right)^{\alpha-\tilde{\delta}-1}$$

$$\geq (1-5\alpha\delta)3^{-\alpha}|\ln\varepsilon| n_\varepsilon \left(\frac{n_\varepsilon}{k}\right)^{\alpha-\tilde{\delta}-1}.$$

In the crucial step from the second to the third line of the inequality chain we use that the relation $\alpha - \tilde{\delta} > 1$ imposes that the function $x \mapsto x^{\alpha-\tilde{\delta}}$ is *convex*, and, therefore, the minimum is taken for the choice $x_i = \frac{(1-\delta)^2 a}{k\varepsilon}$, $0 \leq i \leq k$, due to (3.35). The same conditions imply the following inequalities which are used in the last line of the chain. In fact, for sufficiently small $\varepsilon$, we obtain

$$\text{(5.43)} \qquad f(g_\varepsilon) \geq (1-\delta)\left(\frac{g_\varepsilon}{q_\varepsilon}\right)^\alpha f(q_\varepsilon) \geq (1-\delta)\left(\frac{\delta}{6}\right)^\alpha |\ln \varepsilon|$$

and

$$\text{(5.44)} \qquad \left(\frac{a}{\varepsilon g_\varepsilon n_\varepsilon}\right)^{\alpha-\tilde{\delta}} \geq \left(\frac{2}{\delta}\right)^{\alpha-\tilde{\delta}} \geq (1-\delta)\left(\frac{2}{\delta}\right)^\alpha.$$

Summarizing our conclusions, we may continue the estimate in (5.41) by the inequality

$$\mathbf{P}\left(\sum_{i=1}^k |\varepsilon W_i| \geq (1-\delta)a\right)$$

$$\text{(5.45)} \qquad \leq \beta_\varepsilon^{-k} \exp\left(-(1-5\alpha\delta)3^{-\alpha}|\ln \varepsilon|n_\varepsilon \left(\frac{n_\varepsilon}{k}\right)^{\alpha-\tilde{\delta}-1} + 2n_\varepsilon \ln |\ln \varepsilon|\right)$$

$$\leq \beta_\varepsilon^{-k} \exp\left(-(1-6\alpha\delta)3^{-\alpha}|\ln \varepsilon|n_\varepsilon \left(\frac{n_\varepsilon}{k}\right)^{\alpha-\tilde{\delta}-1}\right),$$

again valid for $\varepsilon$ small enough uniformly over $1 \leq k \leq n_\varepsilon$.

It remains to include the probabilities $\mathbf{P}(N_T = k)$ for $1 \leq k < n_\varepsilon$ into our estimates. For this purpose, we shall estimate $\max_{1 \leq k < n_\varepsilon} \mathbf{P}(N_T = k)\mathbf{P}(\sum_{i=1}^k |\varepsilon W_i| > (1-\delta)a)$. This will be done by looking separately at the cases $k \in A$ and $k \in B$, where

$$\text{(5.46)} \qquad A := \left\{1 \leq k < n_\varepsilon : \left(\frac{n_\varepsilon}{k}\right)^{\alpha-\tilde{\delta}-1} > 3^\alpha\right\},$$

$$B := \left\{1 \leq k < n_\varepsilon : \left(\frac{n_\varepsilon}{k}\right)^{\alpha-\tilde{\delta}-1} \leq 3^\alpha\right\}.$$

The estimate $\mathbf{P}(N_T = k) \leq \mathbf{P}(\tau_1 \leq T)^k \leq (\beta_\varepsilon T)^k$ is valid for any $k \geq 1$. So we obtain for small enough $\varepsilon$

$$\max_{k \in A} \mathbf{P}(N_T = k)\mathbf{P}\left(\sum_{i=1}^k |\varepsilon W_i| > (1-\delta)a\right)$$

$$\text{(5.47)} \qquad \leq T^{n_\varepsilon} \exp(-(1-6\alpha\delta)|\ln \varepsilon|n_\varepsilon)$$

$$\leq \exp\left(-(1-7\alpha\delta)3a\frac{|\ln \varepsilon|}{\varepsilon q_\varepsilon}\right) < D_\varepsilon^a.$$

Finally, we consider $k \in B$. By choice of parameters and definition of $B$, we have $\inf_{k \in B} \ln k \geq (1-\delta)|\ln \varepsilon|$ for sufficiently small $\varepsilon$. Hence, for $\varepsilon$ small, again by means of Stirling's formula,

$$\mathbf{P}(N_T = k) \leq \frac{(\beta_\varepsilon T)^k}{k!} \leq \beta_\varepsilon^k \exp(-k(\ln k - 1 - \ln T))$$
$$\leq \beta_\varepsilon^k \exp(-(1-2\delta)k|\ln \varepsilon|)$$
$$\leq \beta_\varepsilon^k \exp\left(-(1-2\delta)n_\varepsilon |\ln \varepsilon| \frac{k}{n_\varepsilon}\right). \tag{5.48}$$

This combines with our estimate for the rate of sums of big jumps to the inequality

$$\max_{k \in B} \mathbf{P}(N_T = k) \mathbf{P}\left(\sum_{i=1}^k |\varepsilon W_i| > (1-\delta)a\right)$$
$$\leq \exp\left(-(1-6\alpha\delta)|\ln \varepsilon| n_\varepsilon \left[3^{-\alpha}\left(\frac{n_\varepsilon}{k}\right)^{\alpha-\tilde{\delta}-1} + \frac{k}{n_\varepsilon}\right]\right)$$
$$\leq \exp\left(-(1-6\alpha\delta)|\ln \varepsilon| n_\varepsilon \inf_{y>0}[3^{-\alpha} y^{-(\alpha-\tilde{\delta}-1)} + y]\right). \tag{5.49}$$

It is easy to see that $\inf_{y>0}[3^{-\rho} y^{-(\rho-1)} + y] = \frac{1}{3}\rho(\rho-1)^{-(1-1/\rho)}$ holds for any $\rho > 1$. The mapping $\rho \mapsto \frac{1}{3}\rho(\rho-1)^{-(1-1/\rho)}$ is continuous on $(1, \infty)$. Thus, $\tilde{\delta}$ can be chosen sufficiently small for the following estimate to hold:

$$\inf_{y>0}[3^{-\alpha} y^{-(\alpha-\tilde{\delta}-1)} + y] \geq 3^{-\tilde{\delta}} \inf_{y>0}[3^{-(\alpha-\tilde{\delta})} y^{-(\alpha-\tilde{\delta}-1)} + y]$$
$$\geq (1-\delta) \inf_{y>0}[3^{-\alpha} y^{-(\alpha-1)} + y] = \frac{1-\delta}{3} d_\alpha. \tag{5.50}$$

So we finally get the inequality

$$\max_{k \in B} \mathbf{P}(N_T = k) \mathbf{P}\left(\sum_{i=1}^k |\varepsilon W_i| > (1-\delta)a\right)$$
$$\leq \exp\left(-(1-7\alpha\delta)\frac{d_\alpha}{3}|\ln \varepsilon| n_\varepsilon\right)$$
$$\leq \exp\left(-(1-8\alpha\delta)a d_\alpha \frac{|\ln \varepsilon|}{\varepsilon q_\varepsilon}\right)$$
$$\leq D_\varepsilon^{(1-8\alpha\delta)a}. \tag{5.51}$$

Now combine this with our estimate on $A$, and take into account that we have $n_\varepsilon$ summands of the two types. But by its choice, the factor $n_\varepsilon$ being a power function of $\varepsilon$ does not change the exponential asymptotics in the limit $\varepsilon \to 0$. This completes the proof of both Lemma 5.1 and Theorem 2.2. □

5.4. *Bounded sub-exponential tails. Proof of Theorem 2.3, lower bound.*
Here we essentially proceed as in Section 5.2. Let $\theta \in (\frac{1}{n}, \frac{1}{n-1}]$ for some $n \geq 1$, and $\vartheta := 1 - (n-1)\theta$. Let $0 < \delta < \frac{1}{3}\alpha(1-\alpha)$ and choose $\delta' \in (0, \delta]$ small enough such that $(1-\delta')(\frac{\vartheta - 4\delta' + 3(\delta')^2}{\vartheta})^{\alpha+\delta'} \geq 1 - \delta$.

Consider the decomposition of the process $L^{\varepsilon,\theta}$ into a sum of small and big jump parts with the threshold $g_\varepsilon$ defined in (5.14). Accordingly define

(5.52) $$f_{\varepsilon,\theta}(u) := -\ln \nu_{\varepsilon,\theta}([u, \infty)).$$

Clearly, $f_{\varepsilon,\theta}(u) \geq f(u)$ for $0 < u \leq \frac{\theta}{\varepsilon}$, and $f_{\varepsilon,\theta}(u) = +\infty$ otherwise.

With $k = k_\varepsilon$ from (5.14), we use the covering (5.13) and the estimate (5.15). One can easily see that only the estimate of the probability of the event $\{\sum_{i=1}^k |\varepsilon W_i| \geq 1 - 3\delta'\}$ should be revisited. First we note that this set is empty for $k < n$ and $\delta$ small enough. For $\varepsilon$ small enough, we get (under convention that infimum over the empty set equals $+\infty$) that

(5.53)
$$\inf\left\{\sum_{i=1}^k f_{\varepsilon,\theta}(x_i) : \sum_{i=1}^k x_i = (1-3\delta')(1-\delta')\varepsilon^{-1}, x_i \in [g_\varepsilon, \varepsilon^{-1}]\right\}$$
$$\geq \inf\left\{\sum_{i=1}^k f_{\varepsilon,\theta}(x_i) : \sum_{i=1}^k x_i = (1-3\delta')(1-\delta')\varepsilon^{-1}, x_i \in [g_\varepsilon, \theta\varepsilon^{-1}]\right\}$$
$$\geq \inf\left\{\sum_{i=1}^k f(x_i) : \sum_{i=1}^k x_i = (1-3\delta')(1-\delta')\varepsilon^{-1}, x_i \in [g_\varepsilon, \theta\varepsilon^{-1}]\right\}.$$

These inequalities are trivial for $\theta \geq 1$. Otherwise we note that $f_{\varepsilon,\theta}(y) = +\infty$ for $y > \theta/\varepsilon$. Further, using Potter's bound $f(x) \geq (1-\delta')(\frac{\theta\varepsilon^{-1}}{x})^{\alpha+\delta'} f(\theta\varepsilon^{-1})$ for $g_\varepsilon \leq x \leq \theta\varepsilon^{-1}$, the estimate (3.33) and the inequality

(5.54)
$$\left(\frac{\vartheta - 4\delta' + 3(\delta')^2}{\theta}\right)^{\alpha+\delta'} \frac{f(\theta/\varepsilon)}{f(\vartheta/\varepsilon)} f(\vartheta/\varepsilon)$$
$$\geq (1-\delta')\left(\frac{\vartheta - 4\delta' + 3(\delta')^2}{\vartheta}\right)^{\alpha+\delta'} f(\vartheta/\varepsilon)$$
$$\geq (1-\delta)f(\vartheta/\varepsilon)$$

and the definition of regularly varying functions, we continue the chain of estimate to obtain that the last expression in (5.53) is bigger than

(5.55)
$$(1-\delta')f(\theta\varepsilon^{-1})(\theta^{-1}\varepsilon)^{\alpha+\delta'}$$
$$\times \inf\left\{\sum_{i=1}^k x_i^{\alpha+\delta'} : \sum_{i=1}^k x_i = (1-3\delta)(1-\delta')\varepsilon^{-1}, x_i \in [g_\varepsilon, \theta\varepsilon^{-1}]\right\}$$
$$\geq (1-\delta')f(\theta\varepsilon^{-1})(\theta^{-1}\varepsilon)^{\alpha+\delta'}((n-1)(\theta\varepsilon^{-1})^{\alpha+\delta'}$$

$$+ ((\vartheta - 4\delta' + 3(\delta')^2)\varepsilon^{-1})^{\alpha+\delta'})$$
$$\geq (1-\delta)((n-1)f(\theta\varepsilon^{-1}) + f(\vartheta\varepsilon^{-1})) \geq (1-2\delta)\phi(\theta)f(\varepsilon^{-1}).$$

Thus, the estimate similar to (5.20) of the previous section is established, and the statement of the theorem follows.

**Acknowledgments.** The authors thank the Associate Editor and two anonymous referees for their careful reading of the manuscript and valuable suggestions.

P. IMKELLER
I. PAVLYUKEVICH
T. WETZEL
HUMBOLDT—UNIVERSITÄT ZU BERLIN
DEPARTMENT OF MATHEMATICS
RUDOWER CHAUSSEE 25
12489 BERLIN
GERMANY
E-MAIL: imkeller@math.hu-berlin.de
       pavljuke@math.hu-berlin.de
       wetzel@math.hu-berlin.de